\documentclass[11 pt, reqno]{article}

\usepackage{amsmath,amsfonts,amssymb,amsthm}
\usepackage[colorlinks=true,hyperindex=true]{hyperref}
\usepackage{cancel,bbm}
\usepackage{mathabx} 
\usepackage{comment}
\usepackage{enumitem}
\usepackage{cite}

\usepackage{chngcntr}
\counterwithin*{equation}{section}

\usepackage[a4paper, left=2cm, right=2cm, top=2 cm, bottom= 2 cm]{geometry}
\usepackage{color}
\usepackage{fancyhdr}
\usepackage{latexsym}

\newtheorem{ccounter}{ccounter}[section]
\newtheorem{thm}[ccounter]{Theorem}
\newtheorem{lem}[ccounter]{Lemma}
\newtheorem{cor}[ccounter]{Corollary}
\newtheorem{defn}[ccounter]{Definition}
\newtheorem{prop}[ccounter]{Proposition}
\newtheorem{ass}[ccounter]{Assumption}
\newtheorem{ex}[ccounter]{Example}

\def\bet{\begin{thm}}
\def\eet{\end{thm}}
\def\bel{\begin{lem}}
\def\eel{\end{lem}}
\def\bas{\begin{ass}}
\def\eas{\end{ass}}
\def\bec{\begin{cor}}
\def\eec{\end{cor}}
\def\bed{\begin{defn}}
\def\eed{\end{defn}}
\def\bep{\begin{prop}}
\def\eep{\end{prop}}
\def\beq{\begin{equation}}
\def\eeq{\end{equation}}
\def\proof{\noindent {\bf Proof.}\ \ }
\def\bea{\begin{equation*}}
\def\eea{\end{equation*}}
\def\tr{\mathrm{tr}}
\def\bex{\begin{ex}}
\def\eex{\end{ex}}

\def\F{\mathcal{F}}
\def\rr{\mathbb{R}}

\def\cc{\mathbb{C}}
\def\zz{\mathbb{Z}}

\def\i{\mathrm{i}}
\def\e{\rm{e}}

\def\del{\partial}
\def\d{\mathrm{d}}
\def\e{\mathrm{e}}

\def\Re{\mathrm{Re}}
\def\supp{\mathrm{supp}\,}
\def\Im{\mathrm{Im}\mbox{ }}
\def\pp{\mathbb{P}}
\def\ee{\mathbb{E}}
\def\eps{\varepsilon}
\def\mfctN{m_{\mathrm{fc}, t}^{(N)}}
\def\rhofctN{\rho_{\mathrm{fc}, t}^{(N)}}

\renewcommand\leq\varleq
\renewcommand\geq\vargeq
\def\D{\mathcal{D}}

\def\z1{\zeta_1}
\def\z2{\zeta_2}

\def\z{\zeta}

\def\O{\mathcal{O}}

\def\tt{{\mathbb{T}}}

\def\btt{{(\mathbb{T})}}

\def\btti{{(\mathbb{T}i)}}
\def\bttj{{(\mathbb{T}j)}}
\def\bttk{{(\mathbb{T}k)}}

\def\P{\mathcal{P}}

\def\A{\mathcal{A}}
\def\L{\mathcal{L}}
\def\1{\boldsymbol{1}}

\def\gamjl{\gamma_{j, T}}

\def\Gddi{\mathcal{G}_{\delta, i}}

\def\dist{\mathrm{dist}}

\def\I{\mathcal{I}}
\def\Delone{\Delta_1}

\def\L{L}

\def\Ht{H_t}

\def\remark{\noindent{\bf Remark. }}

\def\rhofcll{\rho_{\mathrm{fc}, T}}
\def\hatH{\hat{H}}
\def\rhofc{\rho_{\mathrm{fc}}}
\def\hatx{\hat{x}}
\def\mfcll{m_{\mathrm{fc}, T}}
\def\haty{\hat{y}}
\def\tilx{\tilde{x}}
\def\tily{\tilde{y}}

\def\bG{\bar{G}}

\def\L{\mathcal{L}}

\def\onegdd{\1_{\mathcal{G}_{\delta, i}}}
\def\onegddc{\1_{\mathcal{G}^c_{\delta, i}}}

\def\Ieps{\mathcal{I}_{\varepsilon}}
\def\G{\mathcal{G}}

\def\lami{\lambda^{(i)}}

\def\Ni{N^{(i)}}
\def\Nj{N^{(j)}}

\def\mfc{m_{\mathrm{fc}}}
\def\mN{m_N}

\def\Lam{\Lambda}
\def\Om{\Omega}

\def\NL{N_L}

\def\UB{\mathcal{U}^{(B)}}
\def\fa{\mathfrak{a}}
\def\rhon{\rho^{(n)}}
\def\mfct{m_{\mathrm{fc}, t}}
\def\rhofct{\rho_{\mathrm{fc}, t}}
\def\om{\omega}

\def\IV{\mathcal{I}_V}

\def\Tom{\mathcal{T}_\omega}
\def\nn{\mathbb{N}}
\def\Fk{(F_{i_k})_{\sigma_k}}
\def\bttij{(\mathbb{T}ij)}

\def\IV{\mathcal{I}_V}
\def\DLmu{\mathcal{D}_{L, \mu}}
\def\mNt{m_{N, t} }
\def\nfct{n_{\mathrm{fc}, t} }
\def\nt{n_t}
\def\gamit{\gamma_{i, t}}

\def\Amut{\mathcal{A}_{q, t}}
\def\Amuto{\mathcal{A}_{q, t_0}}
\def\Ieta{\mathcal{I}_\eta}
\def\tilH{\tilde{H}}

\def\lamit{\lambda_{i, t}}
\def\tilrho{\tilde{\rho}}
\def\tileta{\tilde{\eta}}

\def\Gdelt{\mathcal{G}_\delta}
\def\NL{N_L}

\def\tilW{\tilde{W}}
\def\rhosc{\rho_{\mathrm{sc}}}
\def\rhoscab{\rho_{\mathrm{sc}}^{(a, b)}}
\def\muab{\mu^{(a, b)}}
\def\gamkt{\gamma_{k, t}}
\def\nscab{n_{\mathrm{sc}}^{(a, b)}}
\def\gamjt{\gamma_{j, t}}
\def\gamli{\gamma_{T, i} }

\def\mfcllN{m_{\mathrm{fc}, T}^{(N)}}
\def\rhofcllN{\rho_{\mathrm{fc}, T}^{(N)}}

\def\lamipt{\lambda_{i+1, t}}

\def\mfcl{m_{\mathrm{fc} } }

\def\oneni{\1_{ \{ N^{(i)}_{\mathcal{I}_{\eps} } \geq 1 \}} }
\def\rhofcl{\rho_{\mathrm{fc}}}
\def\mV{m_V}
\def\Iketa{\I_{\eta}^{(k)}}

\def\rhofcto{\rho_{\mathrm{fc}, t_0}}
\def\gamkoto{\gamma_{k_0, t_0} }
\def\hatlam{\hat{\lambda}}
\def\hatmu{\hat{\mu}}
\def\bhatlam{\boldsymbol{\hat{\lambda}}}
\def\blam{\boldsymbol{\lambda}}
\def\bhatmu{\boldsymbol{\hat{\mu}}}
\def\lamil{\lambda_{i, T}}
\def\lamipl{\lambda_{i+1, T}}
\def\pW{p_W}
\def\pG{p_G}
\def\pWn{p_W^{(n)}}
\def\pGn{p_G^{(n)}}

\def\L{G}
\def\sqT{\sqrt{T}}
\def\Ione{\mathcal{I}^{(1)}}
\def\Itwo{\mathcal{I}^{(2)}}
\def\Qone{\mathcal{Q}_1}
\def\tilrho{\tilde{\rho}}
\def\tilm{\tilde{m}}

\begin{document}
\title{Deformed GOE}


\begin{table}
\centering

\begin{tabular}{c}
\multicolumn{1}{c}{\Large{\bf Convergence of local statistics of Dyson Brownian motion }}\\
\\
\\
\end{tabular}
\begin{tabular}{c c c}
Benjamin Landon & & Horng-Tzer Yau\\
\\
 \multicolumn{3}{c}{ \small{Department of Mathematics} } \\
 \multicolumn{3}{c}{ \small{Harvard University} } \\
\small{landon@math.harvard.edu} &  & \small{htyau@math.harvard.edu}  \\
\\
\end{tabular}
\\
\begin{tabular}{c}
\multicolumn{1}{c}{\today}\\
\\
\end{tabular}

\begin{tabular}{p{15 cm}}
\small{{\bf Abstract:} We analyze the rate of convergence of the local statistics of Dyson Brownian motion to the GOE/GUE for short times $t=o(1)$ with deterministic initial data $V$.  
Our main result states that if the density of states of $V$ is bounded both above and away from $0$ down to scales $\ell \ll t$ in a small interval of size $G \gg t$ around an energy $E_0$, then the local statistics coincide with the GOE/GUE near the energy $E_0$ after time $t$.  Our methods are partly based on the idea of coupling two Dyson Brownian motions from \cite{homog}, the parabolic regularity result of \cite{gap}, and the eigenvalue rigidity results of \cite{Kevin1}.}
\end{tabular}
\end{table}

\section{Introduction} \label{sec:int}

{\let\thefootnote\relax\footnote{The work of B.L. is partially supported by NSERC.  The work of H.-T. Y. is partially supported by NSF Grant DMS-1307444 and a Simons Investigators fellowship.}}Wigner ensembles consist of $N \times N$ real symmetric or complex Hermitian matrices $W$ whose entries are random variables that are independent up to the symmetry constraint $W = W^*$.  Wigner's global semicircle law \cite{wigner} states that in the appropriate scaling the empirical distribution of the eigenvalues $(\lambda_i)$ converges to
\beq
\frac{1}{N} \sum_{i=1}^N \delta_{\lambda_i } (E) \to \rhosc (E) := \frac{1}{2 \pi } \1_{ \{ |E| \leq 2 \} }\sqrt{  4 - E^2 }, \quad (E \in \rr )
\eeq
in the weak sense as $N \to \infty$.  The distribution $\rhosc (E)$ is referred to as the semicircle law.  Wigner obtained his global semicircle law by computing the moments $\ee [\tr (W^n)]$ for each $n$.

We denote by $\pW ( \lambda_1, ..., \lambda_N )$ the joint probability density of the unordered eigenvalues of $W$.  If the entries of $W$ are independent real or complex Gaussian random variables with variance equal to $N^{-1}$ then the joint density is explicitly computable and is given by
\beq \label{eqn:pg}
\pG ( \lambda_1, ..., \lambda_N ) = \frac{1}{ Z_G} \prod_{i < j } | \lambda_i - \lambda_j |^\beta \e^{ - \beta N \sum_{i=1}^N \lambda_i^2 /4 }
\eeq
where $\beta$ is $1$ or $2$ for the real and complex cases, respectively.  Above, $Z_G$ is a normalization constant which can be computed explicitly.  These special cases are known as the Gaussian Orthogonal and Gaussian Unitary ensembles (GOE and GUE). The $n$-point correlation functions are defined by
\beq
\pWn  ( \lambda_1, ..., \lambda_n ) := \int_{ \rr^{N-n} } \pW ( \lambda_1, ..., \lambda_N ) \d \lambda_1 ... \d \lambda_N.
\eeq
For the GOE and GUE, the $n$-point correlation functions have been computed explicitly by Dyson, Gaudin and Mehta (see, for example, \cite{mehta2004random}) using orthogonal polynomial techniques exploiting the Vandermonde determinant structure.  

For the case $\beta =2$, the work of Dyson, Gaudin and Mehta asserts that at every fixed energy $E \in (-2, 2)$ in the bulk of the spectrum, 
\beq
\frac{1}{ \rhosc (E)^n } \pGn \left( E + \frac{ \alpha_1}{ N \rhosc (E) } , ...., E + \frac{ \alpha_n }{ N \rhosc (E) } \right) \to \det ( K ( \alpha_i - \alpha_j ) )_{i, j=1}^n
\eeq
where $K$ is the sine kernel
\beq
K (x - y) = \frac{ \sin \pi (x-y) } { \pi (x-y) }.
\eeq
The rescaling by a factor of $N^{-1}$ corresponds to the typical distance between consecutive eigenvalues, and we refer to laws under such a scaling as local statistics.  There are similar but more complicated formulas for the GOE.  

The Wigner-Dyson-Gaudin-Mehta conjecture, or the `bulk universality' conjecture, states that the local eigenvalue statistics of Wigner matrices are universal in the sense that they depend only on the symmetry class of the random matrix ensemble (i.e., real symmetric or complex Hermitian) but are otherwise independent of the underlying law of the matrix entries.  This conjecture has been established for all symmetry classes in the works \cite{erdos2010bulk, erdos2012spectral, erdos2011universality, erdos2012rigidity, general, gap}.  Parallel results were obtained independently in various cases in \cite{tao2011random, tao2010random}.

We remark that when we refer to `bulk universality' in this paper, we refer to the vague convergence of the correlation functions $p_W^{(n)}$ in the averaged energy sense of \cite{erdos2010bulk, erdos2012spectral, erdos2011universality, erdos2012rigidity, general, gap}.  Results at fixed energy for the real symmetric case were recently obtained in \cite{homog} but we will not address this type of convergence in this work.

In order to place the current work in context, we recall the three-step strategy of the proof of bulk universality of Wigner matrices.
\begin{enumerate}[label={(\arabic*)}]
\item \label{item:local} Establish a local semicircle law controlling the density of eigenvalues down to the optimal scale.
\item \label{item:ergod} Prove universality of Wigner matrices with a small Gaussian component by analyzing the convergence of Dyson Brownian motion to local equilibrium.
\item \label{item:green} Prove universality of a general Wigner matrix by comparing its local statistics to an approximating Wigner ensemble with a small Gaussian component.
\end{enumerate}

For an overview of this strategy and a survey of recent results we refer the reader to \cite{localspectral}.  In this paper we are mainly interested in Step \ref{item:ergod}.  The local ergodicity of Dyson Brownian motion (DBM) is the intrinsic mechanism behind the universality of local statistics.

In the present work we analyze the speed of convergence of Dyson Brownian motion for classical values of $\beta$ with deterministic initial data $V$.  Our main result is that if the density of states of $V$ is bounded above and below at all scales down to $\ell \ll t$ in a window of size $G \gg \sqrt{t}$ around an energy $E_0$ then the local statistics of Dyson Brownian motion exhibit bulk universality at time $t$ near $E_0$. We allow for scales as small as $\ell = N^{-1}$, and $G$ can shrink as long as it satisfies $G \gg \sqrt{t}$.

A completely analytic approach to analyzing the time to convergence of Dyson Brownian motion was initiated in \cite{erdos2011universality}, and was further developed in \cite{localrelaxation, erdos2012rigidity, erdos2010universality}.  In these works the optimal rate $t \gtrsim N^{-1}$ was obtained when the initial data is a Wigner matrix.  The key idea was to estimate the entropy flow of Dyson Brownian motion with respect to a global instantaneuous equilibrium state constructed from the invariant semicircle law of the GOE/GUE. 

  For deterministic initial data, the study of the convergence of the local statistics of Dyson Brownian motion was initiated in \cite{Kevin2}.  Under some weak global conditions on the initial data, it was shown that the local statistics coincide with the GOE/GUE for times of order $1$. In this case, the global statistics of Dyson Brownian motion are not close to the semicircle law and are in fact time dependent.  Instead of comparing DBM to a global equilibrium state constructed from the GOE/GUE, an equilibrium state was constructed from a time dependent reference $\beta$-ensemble specifically chosen to match the global eigenvalue density of the DBM.  The analysis of entropy flow with respect to time dependent local equilibrium states was initiated in the work \cite{yau1991relative}.  This method allowed for the comparison of the local statistics of the DBM to that of a $\beta$-ensemble.  Bulk universality for such $\beta$-ensembles was achieved in the series of works \cite{bourgade2012bulk, bourgade2014edge, bourgade2014universality} and therefore the local statistics for DBM was also obtained.  Recently, alternative approaches to the local statistics of $\beta$-ensembles have been presented in \cite{shcherbina} and \cite{transport}.  However, the key input in the proof of \cite{Kevin2} was the universality of the local equilibrium measures which was only proved in \cite{bourgade2012bulk, bourgade2014edge, bourgade2014universality}.

As previously stated, we study the rate of convergence of DBM for short times with deterministic inital data satisfying only a local regularity condition.  It is therefore not possible to compare the DBM to a global equilibrium state as in the works \cite{erdos2011universality, localrelaxation, erdos2012rigidity, erdos2010universality, Kevin2}.  One may attempt to circumvent this and assume regularity globally.  For the class of initial data we are interested in, the local density of the reference $\beta$-ensemble could have order $1$ fluctuations on quite small scales and would therefore be quite rough. In order to complete the approach of \cite{Kevin2}, the  results of \cite{bourgade2012bulk, bourgade2014edge, bourgade2014universality} would have to be extended to $\beta$-ensembles with order $1$ fluctuations over small scales.  While this approach may be possible, it is more appealing to approximate DBM locally.

 We approximate DBM locally by constructing a rescaled and shifted GOE matrix $\tilH$ whose eigenvalues match the DBM in a small window near $E_0$.  We then apply two ideas from previous works.  The first is an idea of \cite{homog}:  we couple the evolution of the eigenvalues of the two ensembles under the Dyson Brownian motion, so that their differences satisfy a system of difference equations which may be interpreted as a random walk in a random environment.  The second idea is the parabolic regularity result for such systems of difference equations of \cite{gap}.  
In order to use the approach outlined here,  two main  ingredients are needed. 
These are a rigidity result for DBM with deterministic initial data after a short time (i.e., the analogue of Step \ref{item:local} above) and a level repulsion estimate.  {\it One of our key observations is that the regularity of the initial data down to scales $\ell$ guarantees that both rigidity and level repulsion occur after time $t \gg \ell$.}

For times of order $1$, rigidity was established in the paper \cite{Kevin1} and later refined in \cite{Kevin2}.  The adaptation to short times is a minor modification of the proofs there, and we will only state which changes are necessary in lieu of a complete proof.

Previously, level repulsion estimates for Wigner ensembles whose entries have a smooth distribution were obtained in \cite{wegner}. Weaker level repulsion estimates but with no smoothness condition were obtained in \cite{tao2011random, taolevel}.  The proof of \cite{wegner} was modified in \cite{homog} to include the case when the Wigner ensemble is not smooth but instead is the sum of a (possibly non-smooth) Wigner matrix and an independent Gaussian part.  However, these estimates degenerate as the Gaussian component becomes small, and would therefore only be useful for our purposes if we were interested in times of the order $t = N^{-\eps}$.  We establish new level repulsion estimates which show that as long as $t \gg \ell$ then one already has (essentially) as much level repulsion as one would have for times of order $1$.   

Putting these ingredients together, we will  prove that at a short time after the coupling is initiated, the eigenvalue gaps of the two ensembles coincide down to a scale $N^{-1-\eps}$ with high probability.  This proves the fixed label gap universality for DBM and also the bulk universality of the $n$-point correlation functions in the aforementioned locally averaged sense.

As an application of our work we prove bulk universality for deformed Wigner ensembles with a small Wigner component. The case of a large Wigner component was proved in \cite{Kevin2}.  Another application of our result is the bulk universality of sparse Erd\H{o}s-R\'enyi graphs.  In a paper with J. Huang \cite{huang} we prove that Erd\H{o}s-R\'enyi graphs where the probability $p$ of each edge occuring is as small as $p \geq N^{\eps}/N$ exhibit bulk universality.  The previous result obtained in \cite{erdos2012spectral, renyione} allowed for $p$ only as small as $p \geq N^{2/3+\eps}/N$.

We outline the rest of the paper.  In Section \ref{sec:results} we define our model and state our main results on bulk and gap universality.  In Section \ref{sec:local} we state the local law for deformed Wigner ensembles and state the rigidity estimates for the eigenvalues.  We also state our level repulsion estimates.  Section \ref{sect:dbm} contains the main novelty of this paper, our analysis of the Dyson Brownian motion with initial data $V$.  In Section \ref{sec:lr} prove our level repulsion estimates.  In Section \ref{sect:mr} we give the proofs of our main results.  In Section \ref{sect:llproof} we state and prove some deterministic facts required for the proof of the local law, and then give a proof of the local law.  We also derive the rigidity estimates from the local law.

After completing this manuscript we learned that similar results were obtained independently in \cite{ES}.

\noindent{\bf Acknowledgements.} The authors thank Roland Bauerschmidt for helpful comments on a preliminary draft of this manuscript.  The authors thank J. Huang for pointing out an improvement in the assumptions of Definition \ref{def:v}.


\section{Main results}\label{sec:results}

Before defining our model and stating our results, we remark that we will only state and prove our results in the real symmetric, i.e., $\beta =1$ case. The adaptation to the complex Hermitian case, i.e., $\beta =2$, requires only notational changes.

\subsection{Definition of model}
In this section we introduce the model under consideration.  

\bed \label{def:v}Let $\ell = \ell_N$ and $\L = \L_N$ be two $N$-dependent parameters satisfying$^1$ 
{\let\thefootnote\relax\footnote{1.  In a previous draft our assumption was $G^2 \gg \ell$ instead of $G \gg \ell$.  The authors thank J. Huang for pointing out this improvement.}}
\beq
\frac{1}{N} \leq \ell \leq N^{-\eps_1}, \qquad N^{\eps_1} \ell \leq \L \leq N^{- \eps_1} \label{eqn:eps1}
\eeq
for some $\eps_1 > 0$.
A deterministic diagonal matrix $V = \mathrm{diag}(V_1, ..., V_N)$ is called $(\ell, \L)$-regular at $E_0$ if there are constants $c_V > 0, C_V > 0$ so that the following holds.  On the interval
\beq
\I_{E_0 , \L} := (E_0 - \L, E_0 + \L )
\eeq
we have
\beq \label{eqn:vass}
c_V \leq \Im \mV (E + \i \eta ) \leq C_V
\eeq
uniformly for $E \in \I_{E_0, G}$ and $\ell \leq \eta \leq 10$, where
\beq
\mV (z) := \frac{1}{N} \sum_{i=1}^N \frac{1}{ V_i - z }
\eeq
is the Stieltjes transform of $V$.  Moreover, we assume that there is a fixed number $B_V >0$ so that
\beq
 |V_i| \leq N^{B_V}
 \eeq
  for every $i$.

\eed

We will use the notation
\beq
\I_{E, G} := ( E - \L, E + \L).
\eeq

\bed The Gaussian Orthogonal Ensemble (GOE) consists of symmetric matrices $W$ whose entries are independent centered Gaussians (up to the constraint $w_{ij} = w_{ji}$) with variance
\beq
\ee [ w_{ij}^2 ] = \frac{1 +\delta_{ij}}{N}.
\eeq

\eed

\bed The deformed GOE consists of symmetric matrices
\beq
H_T := V + \sqrt{T} W
\eeq
where $V$ is a deterministic diagonal matrix and $W$ is a GOE matrix and $T \geq 0$ is a real parameter.

\eed

\remark Up to a trivial constant rescaling which goes to $1$ in the $N \to \infty$ limit, the solution of DBM with initial data $V$ is equal to $H_T$ at time $t$ in law with $T \asymp t$, for $ t  = o(1)$.

\subsection{Semicircle and deformed semicircle laws}

In order to state our results on the local eigenvalue statistics we must introduce the macroscopic eigenvalue densities of the GOE and deformed GOE.  The macroscopic eigenvalue density of the GOE is given by the semicircle law:
\beq
\rhosc (E) := \1_{ \{ |E| \leq 2 \}} \frac{ \sqrt{ 4 - E^2 }}{2 \pi}.
\eeq
While the macroscopic eigenvalue density of Wigner (and generalized Wigner) matrices also follows the semicircle law, the deformed GOE follows a deformation of $\rhosc$, the so-called \emph{free convolution} of the semicircle law and the empirical measure of $V$.  We define it through its Stieltjes transform.  We let $\mfcllN$ be the solution to
\beq
\mfcllN (z) = \frac{1}{N} \sum_{i=1}^N \frac{1}{ V_i - z - T \mfcllN (z) }.
\eeq
Above, $V_i$ denote the entries of $V$.  The properties of the above equation are well-studied.  It is known that there is a unique solution to the above equation and that $\mfcllN$ is the Stieltjes transform of a measure which has a density $\rhofcllN$.  This density is compactly supported and analytic on the interior of its support.  We refer the reader to, for example, \cite{biane} for details.  We emphasize that $\mfcllN$ and $\rhofcllN$ and their qualitative properties depend on $N$.  For example, the density $\rhofcllN$ can become quite rough as $N \to \infty$.  For notational convenience we suppress the superscript and denote
\beq
\mfcll (z) := \mfcllN (z) , \qquad \rhofcll (E) := \rhofcllN (E).
\eeq
When we wish to emphasize the $N$-dependence we will use $\mfcllN$ and $\rhofcllN$ instead.  

To state our result on gap universality we define the classical eigenvalue locations of the GOE and deformed GOE by
\beq \label{eqn:classdef}
\frac{i}{N} = \int_{- \infty}^{\gamli} \rhofcll (E) \d E , \qquad \frac{i}{N} = \int_{-\infty}^{\mu_i } \rhosc (E) \d E.
\eeq
The classical eigenvalue locations of the deformed law depend on $N$ but we again suppress this in our notation.

\subsection{Bulk and gap universality}

For a deformed GOE matrix $H = V + \sqT W$ we denote the $n$-point correlation function by $\rhon_T$.  It is defined by
\beq
\rhon_T ( \lambda_1, ..., \lambda_n ) = \int_{\rr^{N-n} } \rho^{(N)}_{T} ( \lambda_1, ...., \lambda_N ) \d \lambda_{n+1} ... \d \lambda_N
\eeq
where $\rho^{(N)}_{T}$ is the joint density of the unordered eigenvalues of $H_T$.  The $n$-point correlation functions of the GOE are denoted by $\rhon_{GOE}$ and are defined similarly.  For the joint eigenvalue density of the GOE we have the explicit expression

\beq
\rho^{(N)}_{GOE} ( \lambda_1, ..., \lambda_N ) = \frac{1}{Z^{(N)}_{GOE}} \prod_{i < j } | \lambda_i - \lambda_j | \e^{ - N \sum_{i=1}^N \lambda_i^2 /4 }
\eeq
where $Z^{(N)}_{GOE}$ is a normalization constant.  Our main result on bulk universality is the following.

\bet \label{thm:bulk} Let
\beq
H_T = V + \sqT W
\eeq
be a deformed GOE matrix. Suppose that $V$ is $(\ell, \L)$-regular at $E$ and that $N^{-\eps} \L \geq T \geq N^{\eps} \ell$ for some $\eps >0$.  Let $O \in C_0^\infty ( \rr^n)$ be a test function.  
Fix a parameter $b = N^{c}  / N$ for any $c >0$ satisfying $c < \eps / 2$.  We have,
\begin{align}
\lim_{N \to \infty} \int_{E - b}^{E + b} \int_{\rr^n} & O ( \alpha_1, ..., \alpha_n )  \bigg\{ \frac{1}{ (\rhofcllN (E) )^n} \rhon_T \left(  E' + \frac{  \alpha_1}{N \rhofcllN (E)} , ..., E' + \frac{  \alpha_n}{N \rhofcllN (E)} \right) \notag \\
&-  \frac{1}{ (\rhosc (E'') )^n} \rhon_{GOE}\left(  E'' + \frac{  \alpha_1}{N \rhosc (E'')} , ..., E'' + \frac{  \alpha_n}{N \rhosc (E'')} \right) \bigg\} \d \alpha_1 ... \d \alpha_n \frac{\d E' }{2b} = 0
\end{align}
for any $E'' \in (-2, 2)$. 
\eet

\remark The upper bound on the size of the averaging window can be removed as long as it is contained in $\I_{E, \L /2}$ and one replaces $\rhofcll (E)$ by $\rhofcll (E')$.  This is due to the fact that the macroscopic density $\rhofcll$ varies on the scale $T^{-1}$; see Lemma \ref{lem:rhobd} below.

\remark The scaling factor $\rhofcllN (E)$ satisfies $c \leq \rhofcllN (E) \leq C$ for all large $N$; see Lemma \ref{lem:rhobd}.

For the gap universality we have the following.  

\bet \label{thm:gap} Let
\beq
H_T = V + \sqT W
\eeq
be a deformed GOE matrix.  Suppose that $V$ is $(\ell, \L)$-regular at $E$ and that $\L N^{-\eps } \geq T \geq N^{\eps} \ell$ for some $\eps >0$.  There is a $c_\eps >0$ such that the following holds.  Let $O \in C_0^\infty ( \rr^n )$ be a test function. Let $i $ be an index so that the $i$th classical eigenvalue of $H_T$ (defined in (\ref{eqn:classdef})) lies in $\I_{E, \L/2}$. Let $i_1, ..., i_n \in \nn$ with $i_k \leq N^{c_\eps}$ for each $k$.  Then for $N$ large enough we have
\begin{align}
\bigg|  \ee^{(H_T)}& \left[ O \left(  N \rhofcllN ( \gamli)  (\lambda_i  - \lambda_{i+i_1}) , ..., N \rhofcllN ( \gamli) ( \lambda_i  - \lambda_{i+i_n})  \right) \right] \notag \\
&- \ee^{GOE} \left[ O \left(  N \rhosc (\mu_i ) ( \lambda_i  - \lambda_{i+i_1}) , ...,  N  \rhosc ( \mu_i ) ( \lambda_i  - \lambda_{i+i_n} ) \right) \right]  \bigg| \leq N^{- c_\eps}.
\end{align}
\eet

\subsection{Bulk universality for deformed Wigner ensembles}

In this section we summarize how Theorem \ref{thm:bulk} can be used to prove bulk universality for deformed Wigner ensembles with a small Wigner component.  Deformed Wigner ensembles are of the form
\beq
U_T := V + \sqrt{T} W' \label{eqn:dw}
\eeq
where $V$ is a deterministic diagonal matrix and $W'$ is a Wigner matrix.   In \cite{Kevin2}, bulk universality was proven for such ensembles for $T = 1$ by following the three step strategy outlined in Section \ref{sec:int}.  Using our result one can use the three step strategy to show that if the density of states of $V$ is bounded above and away from $0$ down to scales $\ell \ll T$ in a window of size $G \gg \sqT$ around $E_0$, then $U_T$ exhibits bulk universality in the averaged $n$-point correlation sense at the energy $E_0$. As previously stated, we have shown how to adapt the proof of Step \ref{item:local} in \cite{Kevin1, Kevin2} to short times in this paper.  The adaptation of Step \ref{item:green} to deformed ensembles was also achieved in \cite{Kevin2}, and the adaptation to small times is straightforward. We omit it in this paper.  This completes the proof of bulk universality for such ensembles.

The above proof can clearly be extended to the case when $V$ in (\ref{eqn:dw}) is random and independent of $W'$, as long as the regularity properties are uniformly satisfied with high probability and the empirical Stieltjes transforms of $V$ are close to some deterministic Stieltjes transform.  This is a further simple modification of the arguments of \cite{Kevin2} and we omit the details.

\section{Local deformed semicircle law, rigidity and level repulsion} \label{sec:local}
\subsection{Conventions and notations}  In this subsection we introduce some notation and conventions which will be used throughout the paper.  We introduce two parameters $\xi = \xi_N$ and $\varphi = \varphi_N$.  We let
\beq
a_0 < \xi \leq A_0 \log \log  N , \qquad \varphi = ( \log (N) )^{C_1}
\eeq
for some constants $a_0 > 2$, $A_0 \geq 10$ and $C_1 > 1$.

\bed  We say that an event $\Om$ has $(\xi, \nu)$-high probability, if
\beq \label{eqn:param}
\pp [ \Om^c ] \leq \e^{ - \nu ( \log (N) )^{\xi}} , \qquad \nu > 0
\eeq
for all large $N$.  We say that an event $\Om_1$ holds with $(\xi, \nu)$-high probability on an event $\Om_2$ if
\beq
\pp [ \Om_1^c \cap \Om_2 ] \leq \e^{ - \nu ( \log (N) )^\xi }.
\eeq

\eed

\subsection{Deformed semicircle law} \label{sec:deformed}

For reasons that will be clear in Section \ref{sect:dbm} we modify the definition of $H_T$ slightly and introduce
\beq \label{eqn:ht}
H_t := \e^{-t /2} V + ( 1 - \e^{-t} )^{1/2}W.
\eeq
The solution of DBM with initial data $V$ is equal to $H_t$ in law for each $t$.  Up to a trivial rescaling $H_t$ corresponds to the deformed model $H_T$ with $T \asymp t$, for times $t = o (1)$.  Note that the assumptions (\ref{eqn:vass}) on $V$ are essentially invariant under rescaling $V \to \e^{-t/2} V$ for $t = o (1)$.

The deformed semicircle law satisfies the equation
\beq \label{eqn:mfct}
\mfctN (z) = \frac{1}{N} \sum_{i=1}^N \frac{1}{ e^{-t/2} V_i - z - ( 1 - e^{-t} ) \mfctN (z) }.
\eeq
As above, $\mfctN$ is the Stieltjes transform of a measure that has a density $\rhofctN$ that is compactly supported and analytic on the interior of its support.   For notational convenience we will suppress the superscript and use the notation
\beq
\mfct (z) := \mfctN (z), \qquad \rhofct (E) := \rhofctN (E).
\eeq
When we wish to emphasize the $N$-dependence we will use $\mfctN$ and $\rhofctN$ instead.

We now collect some properties of the density $\rhofct$ and $\mfct$.  We will assume that $V$ is $(\ell, \L)$-regular at $E_0$.  For $\om > 0$ satisfying
\beq
\om < \eps_1/10 \label{eqn:omrest}
 \eeq
 where $\eps_1$ is the constant appearing in (\ref{eqn:eps1}) we define the time domain
\beq
\Tom = \{ t : \ell N^{\om } \leq t \leq N^{-\om} \L \}.
\eeq

We defer the proof of the following lemma to Lemma \ref{lem:rstab} in Section \ref{sect:llproof}.

\bel \label{lem:rhobd} Fix $\om > 0$ and $1> q>0$.  There are constants $c, C$ depending only the constants in (\ref{eqn:vass}) such that
\beq
c \leq \rhofct (E) \leq C, \qquad | \rhofct' (E) | \leq \frac{C}{ t}
\eeq
for all $E \in \I_{E_0, q\L}, t \in \Tom$, and $N$ large enough (where large enough depends on $\om$, $q$ and the constants in (\ref{eqn:vass})).

\eel

\subsection{Local law and rigidity}

 For any $L \geq 12 \xi$ and $0 < q < 1$ we define the spectral domain
\begin{align} \label{eqn:spectdomain}
\D_{L, q} := \D_1 \cup \D_2
\end{align}
where
\begin{align} \label{eqn:donedtwo}
\D_1 &:=\{ z = E + \i \eta : E \in \I_{E_0, q\L} , \varphi^L \leq N \eta \leq 10 N \} \notag \\
\D_2 &:= \{ z = E + \i \eta : |E| \leq N^{ 5 B_V} , 10N \leq N \eta \leq N^{10 B_V +2} \} .
\end{align}
 The domain $\D_1$ is where we prove the local deformed law down to the optimal scale $\eta \gtrsim N^{-1}$, using our assumptions on $V$.  On the domain $\D_2$ we are only concerned with the behaviour of $\mN$ when $z$ is either far from the real axis or far from the spectrum of $H_t$.

We define the empirical Stieltjes transform of $\Ht$ as
\beq
\mNt (z) = \frac{1}{N} \tr \left( \frac{1}{\Ht - z } \right).
\eeq

\bet Let $H_t$ be as in (\ref{eqn:ht}) where $V$ is $\ell$-regular on $\I$. Let \label{thm:ll}
\beq
\xi = \frac{ A_0 + o (1) }{2} \log \log N.
\eeq
Fix $0<q<1$.  There are constants $\nu > 0$ and $c_1$ depending on $\IV$, the constants appearing in (\ref{eqn:vass}), $A_0$ and $C_1$ in (\ref{eqn:param}) and the choice of $q$, so that the following holds with $L \geq 40 \xi$.  For any $z \in \D_{L, q}$ and $t \in \Tom$
\beq
| \mNt (z)  - \mfctN (z) | \leq ( \varphi )^{c_1 \xi } \frac{1}{  N \eta}
\eeq
with $(\xi, \nu)$-high probability for $N$ large enough (here, large enough depends on the choice of $\nu$, $\om$ and $q$).

\eet

The proof follows closely that of the same result in \cite{Kevin1, Kevin2}.  We summarize the changes in Section \ref{sect:llproof}.

From this local law we derive the following result on the rigidity of the eigenvalues of $\Ht$. We define the classical eigenvalue location $\gamit$ by 
\beq
\frac{i}{N} = \int_{-\infty}^{\gamit} \rhofct (u) \d u.
\eeq

For a fixed constant $q >0$ satisfying
\beq
0 < q < 1
\eeq
 and a time $t$  we define the bulk index set by
\beq
\Amut := \{ i : \gamit \in \I_{E_0, q \L} \}.
\eeq
We abbreviate
\beq
\lamit := \lambda_i ( \Ht )
\eeq
where $\lambda_i (A)$ denotes the $i$th largest eigenvalue of a matrix $A$.  The following lemma shows that the definition of $\Amut$ is essentially independent of the chosen time $t$.

\bel \label{lem:bulkindices} Let $0 < q_1 < q_2 < 1$  and $\om > 0$. Then for $N$ large enough, we have for all $s, t \in \Tom$,
\beq
\A_{q_1, s} \subseteq \A_{q_2, t}.
\eeq

\eel
We defer the proof to Section \ref{sect:llproof}.

\bet \label{thm:rig} Fix $1 >q > 0$.  There are constants $\nu > 0$ and $c_4 > 0$ depending on $\IV, A_0, a_0$, $C_1$ and the choice of $q$ such that the following holds.  For $i \in \Amut$ we have
\beq
| \lamit - \gamit | \leq (\varphi)^{ c_4 \xi } \frac{1}{N}
\eeq
with $(\xi, \nu)$-high probability for all $t \in \Tom$.

\eet

Deducing Theorem \ref{thm:rig} from Theorem \ref{thm:ll} is a minor modification of the usual approach given in, for example, \cite{general}, and so we will not provide all of the details.  This is sketched in Section \ref{sect:rigsketch}.

\subsection{Level repulsion estimates}

The following level repulsion estimates will be an important tool for our proof.  The proof is given in Section \ref{sec:lr}.

\bet \label{thm:lr} Let $\delta >0$ and $r>0$ and let $\lamit$ denote the eigenvalues of the deformed GOE matrix $H_t$ with $V$ $(\ell, \L)$-regular at $E_0$.  Fix $\om >0$ and $1> q >0$.  For $t \in \Tom$ and  $i \in \A_{q, t}$ we have:
\beq
\pp [ | \lamit - \lambda_{i+1, t} | \leq \eps/N ] \leq N^{\delta} \eps^{2-r}  \label{eqn:levelr1}
\eeq
\eet

\section{Optimal speed of convergence of Dyson Brownian motion}
\label{sect:dbm}

\subsection{The approximating GOE ensemble and statement of result}

In this section we implement our main strategy in proving gap universality with a fixed label by approximating $H_t$ by a GOE matrix around a fixed eigenvalue.  Let $\tilW$ be a GOE matrix.  We denote the semicircle law by
\beq
\rhosc (x) := \1_{ \{ |x| \leq 2 \}} \frac{ \sqrt{ 4 - x^2 }}{2 \pi}
\eeq
and the semicircle law for $a \tilW + b$ by
\beq
\rhoscab (x):= 
\1_{ \{ |x-b| \leq 2 a \} } \frac { \sqrt{4 a^2 - ( x - b )^2} }{ 2 a^2 \pi }.
\eeq
The classical eigenvalues $\muab_k$ are defined by
\beq
\frac{k}{N} = \int_{-2 a +b }^{\muab_k} \rhoscab (u) \d u
\eeq
and we abbreviate $\mu_k = \mu_k^{(1, 0)}$.  The following lemma shows that at a fixed time we can locally approximate the classical eigenvalues of $H_t$ by those of a GOE matrix.  We omit the proof as it is elementary.  The only ingredients are the bounds of Lemma \ref{lem:rhobd}.

\bel \label{lem:matching} Let $k_0 \in \Amut$ be a bulk index, and $t \in \Tom$.  Fix $\alpha >0$ and let $j_0$ satisfy $\alpha N \leq j_0 \leq (1 - \alpha ) N$.  Let $a$ and $b$ be
\beq
 a = \frac{ \rhosc ( \mu_{j_0} ) }{ \rhofct ( \gamma_{k_0, t} ) } , \qquad   b = \gamma_{k_0, t} - a \mu_{j_0} .
\eeq
Then there is a constant $C$ depending only on $\alpha$ and the constants in (\ref{eqn:vass}) so that
\beq
| \muab_{j + j_0} - \gamma_{j + k_0, t} | \leq \frac{C}{N}, \qquad \mbox{ for } | j | \leq \sqrt{Nt}
\eeq
\eel

In the following we will implement a coupling argument to match the eigenvalue gaps of $H_t$ by those of a GOE matrix down to a scale $N^{ -1-\eps}$.   We define $\hatlam_i$ as the solutions to the SDEs
\beq
\d \hatlam_i (t)= \sqrt{\frac{2}{N} } \d B_i + \left( \frac{1}{N} \sum_{k\neq i} \frac{1}{ \hatlam_i (t) - \hatlam_k (t) } - \frac{1}{2} \hatlam_i (t) \right) \d t \label{eqn:hat1}
\eeq
with initial data
\beq
\hatlam_i (0) = \lambda_i (V)
\eeq
the ordered eigenvalues of $V$.  For each fixed time $t$ the vector $\bhatlam$ is distributed as $\blam (H_t )$.  Fixing a time $t_0 > 0$ we define $\hatmu_i (t)$ for $t \geq t_0$  as the solution to the SDE
\beq
\d \hatmu_i (t)= \sqrt{\frac{2}{N} } \d B_{i +k_0 - j_0} + \left( \frac{1}{N} \sum_{k\neq i} \frac{1}{ \hatmu_i (t) - \hatmu_k (t) } - \frac{1}{2} \hatmu_i (t) \right) \d t \label{eqn:hat2}
\eeq
with initial data $\hatmu_i (t_0 ) = \lambda_i ( a \tilW + b)$ with $\tilW$ a GOE matrix and constants $a$ and $b$ which will be chosen according to Lemma \ref{lem:matching}. The Brownian motions appearing above are the same as in equation (\ref{eqn:hat1}).  This idea of coupling two DBMs originated in \cite{homog} and is crucial to our work.  

For fixed times $t$ the vector $\bhatmu$ is distributed as the eigenvalues of $a_t \tilW +b_t$ with deterministic constants $a_t$ and $b_t$ which are easily computed and satisfy $c \leq a_t \leq C$ and $|b_t | \leq C$.

The goal of the remainder of the section is the proof of the following theorem on the differences of the gaps $ ( \hatlam_i - \hatlam_j ) (t)- ( \hatmu_{i + j_0 - k_0 } - \hatmu_{j +j_0 - k_0} ) (t)$ for a time $t = t_0 + N^{\eps -1}$.

\bet \label{thm:matching} There exists universal constants $\fa > 0$ and $c_1 > 0$ such that the following holds.  Let $V$ be $(\ell, G)$-regular at $E_0$.  Let $\bhatlam (t)$ be the solution to the SDE (\ref{eqn:hat1}) with initial data $V$ as described above.  Fix $\om >0$, $1> q>0$ and let $t_0 = N^{\tau_0} / N \in \Tom$.  Fix a bulk index $k_0 \in \Amuto$, $\alpha >0$ and choose $a, b, j_0$ as in Lemma \ref{lem:matching} for $t = t_0$, and let $\tilW$ be a GOE matrix independent of $H_{t_0}$.  Let $\bhatmu$ be the solution of (\ref{eqn:hat2}) with initial data $a \tilW + b$ at time $t_0$ as described above.

Let $0 < \om' \leq \tau_0$.  Then, there is an event $\F$ with probability $\pp [ \F ] \geq 1 - N^{-c_1 \om'}$ on which the following holds.
\begin{align}
| \hatlam_{j +k_0} (t_0 + t )- &\hatlam_{j' + k_0} (t_0 + t) - ( \hatmu_{j + j_0} ( t_0 + t ) - \hatmu_{j' + j_0} ( t_0 + t) ) | \leq \frac{1}{ N^{ \fa \om' +1 } },\notag \\
& \mbox{ for } |j  | + |j' | \leq N^{ \om' / 4000} , \mbox{ and } t \in \left[ \frac{N^{ \om' / 10 } - N^{ \om' / 30}}{N} , \frac{N^{ \om' / 10 }}{{N}} \right]
\end{align}

\eet

Over the next few subsections we will establish the estimates we will need to prove Theorem \ref{thm:matching}.  In the last subsection we will provide the proof.  The strategy is roughly as follows.  Due to the fact that we have coupled the DBM of both ensembles, their differences satisfy a system of difference equations.  The H\"older regularity result for such parabolic systems of \cite{gap} states that such solutions (neglecting edge effects) $v$ satisfy an estimate of the form $|v_i (t) - v_{i+1} (t) | \leq N^{-\eps} ||v (0) ||_\infty$ after a short time $t$.  If we knew that the classical eigenvalue locations of our two ensembles matched throughout the spectrum then by the optimal rigidity results we would have $||v (0 ) ||_\infty \leq N^{-1+\eps/2}$ with overwhelming probability, and we would deduce that the gaps of the two ensembles coincide.  However, as we only have the matching for roughly $N^c$ eigenvalues we will have to cut off  the non-matching eigenvalues.  While this cut-off disrupts the location of the eigenvalues, the gaps remain unchanged down to a scale $N^{-1-\eps}$.  However, due to the cut-off, the differences of the two eigenvalues will only be an approximate solution to the system of difference equations that we have the parabolic regularity result for.  The error term will be controlled by the finite speed of propagation estimates obtained in \cite{gap}, and we will see that the parabolic regularity result in fact applies to our approximate solution as well.  We therefore obtain the desired statement about the gaps.

\subsection{Rescaling, time shift and regularization}

It is natural to rescale the problem and work in the microscopic scaling.  We remark that we only work in microscopic coordinates for the proof of Theorem \ref{thm:matching}.  In all other sections of this paper we state and prove our results in the macroscopic scaling in which the typical distance between consecutive eigenvalues is $O (1/N)$.  For this proof we also introduce a time shift and a relabeling which sets $t_0 = k_0 = j_0 =0$.  To be more precise, we define the variables
\beq
x_j (t) = N \hatlam_{k_0  + j } (t_0 + N^{-1} t), \qquad y_j (t) = N \hatmu_{j_0 + j } ( t_0 + N^{-1} t ).
\eeq
They satisfy the stochastic differential equations
\begin{align}
\d x_k (t) &= \sqrt{2} \d B_k + \sum_{j \neq k} \frac{1}{ x_k (t)- x_j (t) } \d t - \frac{1}{2N} x_k (t) \d t \\
\d y_k (t) &= \sqrt{2} \d B_k + \sum_{j \neq k} \frac{1}{ y_k (t)- y_j (t) } \d t - \frac{1}{2N} y_k (t) \d t
\end{align}
where $B_k$ are independent Brownian motions. Note that the $B_k$ for the $x_k$'s and $y_k$'s are the same.

We first show the existence of a set  $\Gdelt$ of good paths on which rigidity holds for the eigenvalues $x_k$ and $y_k$.
\bel \label{lem:paths}  Let $\om'$ be as in the statement of Theorem \ref{thm:matching}.  
For any $\delta >0$ there is an event $\Gdelt$ with $\pp[ \Gdelt] \geq 1 - N^{-100}$ for all large $N$ such that the following holds.
\begin{align}
c | i - j |  \leq | x_i (t) - x_j (t) | \leq C | i - j | , \qquad | i - j | \geq N^{\om' \delta}, |i| + |j| \leq G N^{1 - \om/20} , 0 \leq t \leq G N^{ 1 - \om/2 } \\
| x_i (t) - x_j (t) | \leq C N^{\om' \delta} | i - j | , \qquad |i | + | j | \leq G N^{1 - \om/20} , 0 \leq t \leq G N^{1 - \om/2},  
\end{align}
and also the same estimates with the $x_j$'s replaced by $y_j$'s.  The constants $c, C$ depend only on the constants appearing in (\ref{eqn:vass}).  Also on $\Gdelt$
\beq
\Im [ m_N (E + \i \eta ) ] \leq C, \qquad E \in \I_{q, G} , \mbox{ } 5 \geq \eta \geq N^{\delta/2}/N
\eeq
and the same for the Stieltjes transform of the matching GOE ensemble.
\eel

\proof The result for the $x_j$'s for fixed times $t$ on a set of high probability follows from the rigidity estimates of Theorem \ref{thm:rig}.  The extension to all times is a minor stochastic continuity argument almost identical to the proof of Lemma 9.3 in \cite{gap}, with the input of the level repulsion bounds of Theorem \ref{thm:lr}.   The result for the $y_j$'s is similar. \qed

We will also need to regularize the dynamics.
 We define the regularized eigenvalues as the solution to the equations
\beq
\d \hatx_k = \sqrt{2} \d B_k +  \sum_{j \neq k} \frac{1}{ x_k (t)- x_j (t) +\eps_{kj} } \d t - \frac{1}{2N} \hatx_k (t) \d t, \qquad \hatx_k (0) = x_k (0), 
\eeq
The constant $\eps_{kj}$ is defined as $\eps_{kj} = \eps$ for $ k > j$ and $\eps_{kj} = - \eps$ for $k < j$.  Here $\eps  = N^{ - 500}$. 
We make the same definition for the $y_k$'s and introduce the regularized $\haty_k$'s.

The following lemma is based on \cite{homog}.  The main inputs are the level repulsion estimates stated in Theorem \ref{thm:lr}. 
\bel \label{lem:regularized} Choose $\eps$ as 
\beq
\eps = N^{-500  }
\eeq
Then there is an event $\F_1$ with probability $\pp[ \F_1] \geq 1 - N^{-100}$ on which
\beq
| \hatx_k (t) - x_k (t) | \leq N^{-20} 
\eeq
for every $k$ with $k+k_0 \in \Amuto$ 
and all $0 \leq t \leq G^2 N^{1 - \om/2}$.

\eel
\proof   
The rescaled difference $q_i (t) := \e^{ t / (2N)} ( x_i (t) - \hatx_i (t) )$ satisfies the equation
\beq
\d q_i = \Om_i (t) \d t , \mbox{ with } \Om_i (t) = \e^{ t / (2N )} \sum_{j \neq i} \frac{ \eps_{ij}}{ ( x_i (t) - x_j (t) ) ( x_i (t) - x_j (t) + \eps_{ij} ) }.
\eeq
We bound
\begin{align}
\ee \sup_{0 \leq t \leq G^2N^{1 - \om/2 } } \left| \int_0^t \Om_i (s) \d s \right| &\leq N^2 \sup_{0 \leq t \leq G^2 N^{1 - \om/2 } }  \sup_{j \in \{i-1, i\} } \left( \ee \frac{1}{ |x_i - x_{j} |^{3/2} }  \right)^{2/3} \left( \ee \frac{\eps^3}{ |x_i - x_{j} + \eps_{ij} |^3 } \right)^{1/3}.
\end{align}
By (\ref{eqn:levelr1}), as we have restricted our attention to $i +k_0 \in \Amut$, we have
\beq
\left( \ee \frac{1}{ |x_i - x_{i+1} |^{3/2} }  \right)^{2/3} \leq N
\eeq
and for any $r>0$
\beq
 \left( \ee \frac{\eps^3}{ |x_{i+1} - x_i + \eps |^3 } \right)^{1/3} \leq  \left( \ee \frac{\eps^{2-r}}{ |x_{i+1} - x_i + \eps |^{2-r} } \right)^{1/3} \leq N \eps^{(2-r)/3} | \log \eps |^{1/3}.
\eeq
The claim then follows from the Markov inequality and our choice of $\eps$. \qed

For the $y_k$'s we have even better level repulsion and so Lemma \ref{lem:regularized} holds for the $y_k$'s as well. 

\subsection{Cut-off of non-matching eigenvalues}
 
   Let now $\om' > 0$ be as in the statement of Theorem \ref{thm:matching} and let
\beq
K = N^{ \om' /2}.
\eeq
Define $\tilx_j$ for $|j| \leq K$ by
\beq
\d \tilx_j (t)= \sqrt{2} \d B_j + \sum_{l \neq k, | l| \leq K} \frac{1}{ x_j (t)- x_l (t) +\eps_{jl} } \d t - \frac{1}{2N} \tilx_j (t) \d t , \quad \tilx_j (0 ) = x_j (0),
\eeq
and a similar definition for $\tily_j$.  We prove the following lemma to control the effect of the cut-off.  The main inputs are the level repulsion estimates of Theorem \ref{thm:lr} and the rigidity estimates of Lemma \ref{lem:paths}.
\bel \label{lem:ot1}  Let $1/2 > \delta >0$, and $t_1 \leq G N^{ 1 - \om/2}$.  There is an event $\F_{2, \delta}$ with $\pp [ \F_{2, \delta} ] \geq 1 - K^{-  \delta/10}$ on which the following estimates hold.  We have,
\beq
\sup_{ 0 \leq t \leq t_1}| \hatx_a (t) - \hatx_b (t) - ( \tilx_a (t) - \tilx_b (t) ) | \leq \frac{ K^{ \delta  } (t_1 +1) | a - b | }{(K - |a| +1)^{1/2} (K - |b|+1)^{1/2} }
\eeq
for all $|a| \leq |b| \leq K$, and then same the estimates changing $\hatx_j$ for $x_j$, and also the same estimates for the $y_j$'s.  

\eel

\proof By Lemmas \ref{lem:paths} and \ref{lem:regularized} we can assume that $\Gdelt$ and $\F_1$ hold.   Define for $|i| \leq K$,
\beq
w_i (t) = e^{ t / (2N ) } ( \hatx_a (t) - \hatx_b (t) - ( \tilx_a (t) - \tilx_b (t) ))
\eeq
to be the rescaled difference.
We have,
\begin{align}
  w_i(t) &= \int_0^t \e^{ s / ( 2N) } \sum_{ |j| > K} \frac{ x_a (s)  - x_b (s) } { ( x_a (s) - x_j (s) + \eps_{aj} ) ( x_b (s) - x_j (s) +\eps_{bj} ) } \d s \notag \\
  &= \int_0^t \e^{ s / ( 2N) } \sum_{ K+ K^{\delta} \geq |j| > K} \frac{ x_a (s)  - x_b (s) } { ( x_a (s) - x_j (s) + \eps_{aj} ) ( x_b (s) - x_j (s) +\eps_{bj} ) } \d s \notag \\
  &+ \int_0^t \e^{ s / ( 2N) } \sum_{ |j| > K +K^{\delta}} \frac{ x_a (s)  - x_b (s) } { ( x_a (s) - x_j (s) + \eps_{aj} ) ( x_b (s) - x_j (s) +\eps_{bj} ) } \d s \notag \\
  &= A_1 + A_2.
\end{align}
We have split the sum into the terms involving $j$ near the edge $K$ and those with $j$ away from the edge $K$.  We will use rigidity to estimate the terms away from $K$ and the level repulsion estimates for those near $K$.  For $A_2$ we get, using Cauchy-Schwartz,
\beq
|A_2| \leq C K^{\delta} |a-b| \int_0^t \left( \sum_{|j| > K + K^{\delta} } \frac{1}{ |x_j(s) - x_a (s) |^2}\right)^{1/2}\left( \sum_{|j| > K + K^{\delta} } \frac{1}{ |x_j(s) - x_b (s) |^2}\right)^{1/2} \d s.
\eeq
We have the estimate
\beq
 \sum_{|j| > K + K^{\delta} } \frac{1}{ |x_j(s) - x_a (s) |^2} \leq \frac{K^{\delta}}{ K - |a| +1} \Im [ m_N (x_a (s) + \i K^{\delta/2} ( K - |a| +1)/N ) ] \leq \frac{K^{2\delta}}{ K - |a| +1} .
\eeq
Hence,
\beq
|A_2| \leq \frac{  K^{3 \delta} |a-b| (t_1 +1 ) }{ ( K - |a| + 1)^{1/2} (K - |b| +1 )^{1/2} } .
\eeq
We estimate $A_1$ by
\beq
|A_1 | \leq \frac{K^{4\delta}| a - b| }{ K - |a| + 1 } \int_0^{t_1} \frac{1}{ (x_{K} (s) - x_{K +1} (s) - \eps)^2 } +\frac{1}{ (x_{-K} (s) - x_{-K -1} (s) + \eps)^2 } \d s  
\eeq
We have estimated the terms $||j| - K | \leq K^\delta$ using, for $j$ near $K$,
\beq
|x_a (s) - x_{j} (s) + \eps_{aj} |\geq  K^{ -  \delta}  ( K - |a| + 1 )|x_{K} (s) - x_{K +1} (s) - \eps |
\eeq
and also 
\beq
| x_b (s) - x_j (s) + \eps_{bj} |\geq  |x_{K} (s) - x_{K +1} (s) - \eps |
\eeq
and then similar estimates for $j$ near $-K$.  We have
\beq
\ee \left[   \int_0^{t_1} \frac{1}{ (x_{K} (s) - x_{K +1} (s) - \eps)^2 } \right] \leq t_1 K^{\delta}
\eeq
by the level repulsion estimate (\ref{eqn:levelr1}) and the same for $-K$.  By the Markov inequality there is an event of probability at least $1 - K^{-\delta}$ on which it less than $t_1 K^{2 \delta}$. This yields the estimate for the $\tilx_k$'s.  Since the event $\F_1$ of Lemma \ref{lem:regularized} holds we also get the estimate for differences of the $x_k$'s and $\hatx_k$'s.  The same proof works for the $y_k$'s. \qed

Let $u_a (t) = \e^{  t / ( 2N ) }( \tilx_a (t) - \tily_a (t) )$, $|a| \leq K$.  For $u$ we have the equation
\begin{align}
\d u_j (t) &= \e^{t / (2N) } \sum_{l \neq j, |l | \leq K} B_{jl} (t) ( x_j (t) - x_l (t) - ( y_j (t) - y_l (t) ) ) \d t \notag \\
&=  \sum_{l \neq j, |l | \leq K} B_{jl} (t) (u_j (t) - u_l (t) ) \d t + \xi_j (t) \d t
\end{align}
where 
\beq
\xi_j (t) = \e^{t / (2N) } \sum_{l \neq j , |l |\leq K}  B_{jl} (t) ( (x_j (t) - x_l (t) ) - ( \tilx_j (t) - \tilx_l (t) ) -(y_j (t) - y_l (t) ) - ( \tily_j (t) - \tily_l (t)  )
\eeq
and
\beq
B_{jl } (t) = \frac{1}{ (x_j (t) - x_l (t) + \eps_{jl} )  (y_j (t) - y_l (t) + \eps_{jl} ) } .
\eeq
Let $v_j$, $|j| \leq K$ be the solution to 
\beq \label{eqn:ot3}
\d v_j (t) = \sum_{l \neq j, |l| \leq K} B_{jl} (t) (v_j (t) - v_l (t) ) \d t 
\eeq
with initial data
\beq
v_j (0) = u_j (0).
\eeq
By the Duhamel formula we have
\beq \label{eqn:ot9}
u_a (t) - v_a (t) = \int_0^t \sum_{|p| \leq K} \UB_{a p } (s, t) \xi_p (s) \d s
\eeq
where $\UB$ is the semigroup associated to the equation (\ref{eqn:ot3}), and $\UB_{ij}$ denote its matrix elements.  To estimate the RHS we require the following two lemmas.  
\bel \label{lem:ot2} Let $\delta >0$ and $0 \leq t_1 \leq G N^{1 - \om/2}$.  There is an event $\F_{3, \delta}$  with probability at least $\pp [ \F_{3, \delta} ] \geq 1 - K^{- \delta/50}$ so that the following estimates hold.  
\beq
 | \xi_a (t) | \leq \frac{ (1 +t_1) K^{\delta}}{ K - |a| + 1} ( | B_{a, a+1}(t) | + |B_{a, a-1} (t) | )
\eeq
for every $|a| \leq K$ and $0 \leq t \leq t_1$.
\eel

\proof Let $\delta >0$.  We can assume that the events $\Gdelt$, $\F_1$ and $\F_{2, \delta}$ of Lemmas \ref{lem:paths}, \ref{lem:regularized} and \ref{lem:ot1} hold.  Proceeding directly from the definition, we estimate the numerator using Lemma \ref{lem:ot1} and obtain
\begin{align}
| \xi_a (t) | &\leq 2 (1 +  t_1) \frac{ K^{\delta}}{ (K - |a| + 1)^{1/2} } \sum_{j \neq a , |j| \leq K} B_{a, j} \frac{| a - j | }{ (K- |j|+1)^{1/2} }\notag \\
&\leq  2 (1 +  t_1) ( |B_{a, a+1}| + |B_{a, a-1} | ) \bigg( \frac{ K^{\delta}}{ K - |a| + 1 }  K^{2\delta}  \notag \\
&+ \frac{K^{\delta}}{ (K - |a|+1)^{1/2} }\sum_{| j - a | \geq K^{\delta} }\frac{1}{ |a - j | ( K - |j|+1)^{1/2} } \bigg) 
\end{align}
where we have used Lemma \ref{lem:paths} in the  last step to estimate the terms $| a - j | \geq K^{\delta}, |j| \leq K$.  We now estimate the sum
\begin{align}
\sum_{  |j| \leq K }\frac{1}{ |a-j | ( K - |j| +1 )^{1/2} } &\leq \sum_{ |j| \leq K, K - |a| +1 \geq 2 (K - |j| +1 ) } \frac{1}{ |a-j | ( K - |j| +1 )^{1/2} } \notag\\
&+\sum_{ |j| \leq K, K - |a| +1 \leq 2 (K - |j| +1 ) , j \neq a} \frac{1}{ |a-j | ( K - |j| +1 )^{1/2} } =:A_1 + A_2.
\end{align}
Clearly,
\beq
A_2 \leq \frac{C}{ (K - |a| +1 )^{1/2} } \log (K).
\eeq
For $A_1$ note that if $( K - |a| +1 )/2 \geq (K - |j| +1)$ then $|j| - |a| \geq (K - |a| +1)/2$ and so
\begin{align}
A_1 &\leq \frac{C}{ (K- |a|+1)^{1/2}} \sum_{|j| \neq a, |j| \leq K} \frac{1}{ |a-j|^{1/2} (K - |j|+1 )^{1/2} } \notag \\
&\leq \frac{C}{ (K- |a|+1)^{1/2} } \sum_{ |j| \neq a , |j| \leq K } \frac{1}{ |a-j| } + \frac{1}{ K - |j|+1} \leq \frac{ C \log(K) }{ ( K - |a| +1)^{1/2} }.
\end{align}
This yields the claim.
\qed 

\subsection{Finite speed of propogation and H\"older regularity}
We require the following finite speed of propogation estimate, the proof of which is a minor modification of the proof of Lemma 9.6 of \cite{gap}.
\bel \label{lem:ot3} Fix $\delta_1 >0$ and $\delta_2 > 0$ and $0 \leq t_1 \leq G N^{1 - \om/2}$.  Suppose that 
\beq \label{eqn:finite1}
B_{ij} \geq \frac{K^{ - \delta_1}}{ |i - j|^2}
\eeq
and
\beq \label{eqn:finite2}
\sup_{0 \leq M \leq K } \frac{1}{1 + t_1 } \int_0^{t_1} \frac{1}{M} \sum_{ |i| \leq M} \sum_{ | j| \leq M} B_{ij} (s) \d s \leq C K^{ \delta_2}
\eeq
hold.  Then for any $0 \leq s \leq t \leq t_1$ we have the estimate
\beq
\UB_{ap} (s, t ) \leq \frac{ C K^{\delta_2 + 2 \delta_1 + 1/2 } \sqrt{ t_1+ 1 } } { | p -a| +1 }
\eeq
\eel

We can now estimate the difference between $u (t)$ and $v(t)$.  Again, we will rely crucially on the level repulsion estimates of Theorem \ref{thm:lr}.
\bel \label{lem:uvdif} Let $\delta > 0$ and $0 \leq t_1 \leq G N^{1 - \om/2}$.  For every $|a| \leq K^{1/2}$, there is an event $\F_{4, \delta, a}$ with $\pp [ F_{4, \delta, a} ] \geq 1 - K^{ - \delta / 500}$ on which
\beq
\sup_{0 \leq t \leq t_1 } | u_a (t) - v_a (t) | \leq \frac{ (1+ t_1)^{5/2} K^{\delta} }{ K^{1/2} }.
\eeq

\eel

\proof In order to apply Lemma \ref{lem:ot3}, we must check that the conditions (\ref{eqn:finite1})-(\ref{eqn:finite2}) are satisfied on a set of high probability.  Without the supremum,  (\ref{eqn:finite2}) with exponent $\delta >0$ holds on a set of probability at least $1 - K^{ - \delta/3}$ for fixed $M$.  The supremum can be replaced by a dyadic choice of  $M = 2^{k}$ for $k\leq C \log N$.  Hence we have that (\ref{eqn:finite2}) holds with exponent $\delta >0$ on a set of probability at least $1 - K^{- \delta /4}$.  The condition (\ref{eqn:finite1}) is a consequence of rigidity (Lemma \ref{lem:paths}) which holds on a set of much higher probability.  Let $\F_4$ be the event that all of this holds as well as the event $\F_{3, \delta}$ of Lemma \ref{lem:ot2}.  Then $\pp [ \F_4] \geq 1 - K^{ - \delta / 60}$.  From (\ref{eqn:ot9}) we have, for all times $0 \leq t \leq t_1$,
\beq
| u_a (t) - v_a (t) | \leq \sum_{|p| \leq K} \int_0^{t_1} | \UB_{a, p} (s, t) | | \xi_p (s) | \d s. \label{eqn:ot10}
\eeq
Using the estimates that hold on $\F_4$ provided by the previous two lemmas we get
\begin{align}
\ee \left[ \1_{\F_4} \sum_{|p| \leq K } \int_0^{t_1} | \UB_{a, p} (s, t_1) \xi_p (s) |    \right] &\leq 2 ( 1 + t_1 )^{5/2}  K^{4 \delta + 1/2}\sum_{ |p| \leq K+1} \sup_{0 \leq t \leq t_1} \frac{ \ee [ B_{p, p+1} (t) ] }{ ( | a - p | + 1 ) ( K - |p| + 1 ) } \notag \\
&\leq ( 1 + t_1 )^{5/2}  \frac{ K^{ 5 \delta }}{ K^{1/2} }.
\end{align}
Hence, there is an event of probability at least $1 - K^{ - \delta}$ on which the RHS of (\ref{eqn:ot10}) is less than $K^{6 \delta-1/2} ( 1 + t_1)^{5/2}$.  This yields the claim. \qed

We now state the H\"older regularity theorem (Theorem 10.3 of \cite{gap}) from which Theorem \ref{thm:matching} will follow.  In order to do so, we need to introduce the following two conditions.

(C1)$_{\rho_1}$ At time $\sigma$ we have
\beq
\sup_{0 \leq s \leq \sigma' } \sup_{1 \leq M \leq K } \frac{1}{ 1 + | s - \sigma' | } \int_s^{\sigma'} \frac{1}{M} \sum_{ |i| \leq K} \sum_{ |j| \leq K } B_{ij} (u) \d u \leq K^{\rho_1}
\eeq
for every $\sigma' \in \{ \sigma \Xi + \sigma \}$
where
\beq
\Xi = \{ - 2^{-m} ( 1 + 2^{-k} )  :  0 \leq k \leq m \leq C \log N \}
\eeq

(C2)$_{\rho_2}$ For every $ 0 \leq s \leq \sigma$ we have 
\beq
B_{ij} (s) \geq \frac{K^{ - \rho_2}}{ | i - j |^2} 
\eeq
and
\beq
\frac{c}{ | i - j |^2} \leq B_{ij} (s) \leq \frac{C}{ | i - j |^2}, \mbox{ for } | i - j | \geq C' N^{\rho_2}
\eeq

\bet \label{thm:holder} There is a universal constant $\fa >0$ such that the following holds. Let $v(t)$ be a solution to the equation (\ref{eqn:ot3}). Let $c_3 >0$ be given, and let $\sigma \in [K^{c_3}, K^{1-c_3} ]$.  Suppose that $(C1)_{\rho_1}$ and $(C2)_{\rho_2}$ hold for $\rho_1$ and $\rho_2$ small enough depending on $c_3$. Then we have
\beq
\sup_{ |j| + |j'| \leq \sigma^{2/3} } | v_j ( t ) - v_{j'} (t) | \leq C \sigma^{ - \fa } || v (0 ) ||_{\infty}
\eeq
for all $t \in [\sigma - \sigma^{1/3} , \sigma ]$.

\eet

\subsection{Proof of Theorem \ref{thm:matching}}

Recall our choice of $K = N^{ \om'/2}$. Let $t_1 = K^{1/10}$ and $\delta = 1/10$.  By Lemmas \ref{lem:ot1} and \ref{lem:uvdif} there is an event $\F_4$ with probability $\pp [ \F_4 ] \geq 1 - K^{ -1/2000}$ on which
\beq
\sup_{ |j| + | j'| \leq K^{1/2000} } \sup_{0 \leq t \leq t_1 } | x_j (t) - x_{j'} (t) - ( y_j (t) - y_{j'} (t) ) - ( v_j (t) - v_{j' } (t) ) | \leq \frac{1}{ K^{1/10}},
\eeq
where $v_j$ was defined in (\ref{eqn:ot3}).  We now apply Theorem \ref{thm:holder} to the difference $v_j (t) - v_{j'} (t)$.  We take $c_3 = 1/10$.  By Lemma \ref{lem:paths} we know that $(C2)_{\rho_2}$ occurs on a set of high probability.  The other condition $(C1)_{\rho_1}$ occurs on a set of probability at least $1 - K^{ - \rho_1 / 4}$. If there were no double supremum this would follow immediately from the level repulsion bounds (\ref{eqn:levelr1}) and Markov's inequality. The double supremum can be replaced by a dyadic choice of $s = 2^{-a} \sigma$ and $M = 2^b$ for integers $a, b \leq C \log N$. Hence we see that there is a constant $c_1$ s.t. on an event of probability at least $1 - K^{ - c_1}$ we have 
\beq
\sup_{ |j| + | j'| \leq K^{1/30} } | v_j (t) - v_{j'} (t) | \leq C K^{ - \fa / 10 } || v ( 0) ||_\infty.
\eeq
for $t \in [ K^{1/10} - K^{1/30} , K^{1/10} ]$. By the rigidity estimates for $H_{t_0}$ and $\hatH_{t_0}$ and our choice of $a$ and $b$ (recall Lemma \ref{lem:matching}) we have that $|| v (0) ||_\infty \leq K^{ \fa / 20}$ on a set of high probability.  This concludes the proof of the theorem. \qed

\section{Level repulsion estimates: proof of Theorem \ref{thm:lr}} \label{sec:lr}
In this section we prove Theorem \ref{thm:lr}.  Previous level repulsion estimates for Wigner ensembles were obtained in \cite{wegner}.  These estimates required a smoothness condition for the law of the entries of $W$.  The dependence of the estimates on the smoothness of the matrix entries was analyzed in the work \cite{homog}.  As a result, the smoothness condition could be relaxed at the cost of introducing an extra error prefactor in the level repulsion estimates depending on the size of the smoothness component of the matrix entries.

Due to the smoothness condition on the matrix entries, the estimates obtained in \cite{wegner, homog} are not sufficient for our purposes.  If we were interested only in times $t = N^{-\eps}$ for small $\eps >0$, then these estimates would suffice.  As we are interested in much smaller times, we must obtain new estimates. Our overall strategy is similar to the method of \cite{wegner} but we will take advantage of the fact that the Wigner part of our ensemble is a GOE matrix which allows us to explicitly compute the expectations of Green's function elements that arise in the proof.

We remark that a form of level repulsion was obtained in \cite{tao2011random, taolevel} for Wigner ensembles under no smoothness assumptions on the matrix elements.  However, the estimates obtained are not strong enough for our methods.  As an additional remark, we note that gap universality implies a weak form of level repulsion; still, such an estimate resulting from gap universality would not be enough for our proof.

The estimates we will obtain are not completely optimal but will suffice for our purposes.  In this section, it will be more notationally convenient to work with the definition
\beq
H_T = V + \sqT W.
\eeq

Before we proceed we remark that we have switched back to the macroscoping scaling, so that the typical distance between eigenvalues is $O (1 / N)$, instead of the microscoping scaling which was only used in the proof of Theorem \ref{thm:matching}.

Theorem \ref{thm:lr} will be derived as a consequence of the following theorem.  For a $1/2 > \eps >0$ and $E \in \I_{E_0, q G}$ let
\beq
\Ieps := [ E - \eps /N, E + \eps / N]
\eeq
\bet \label{thm:lrest} For $k=1, 2$, and $1 > q > 0$ and all $ \delta> 0, r > 0$ we have for $E \in \I_{E_0 , q G }$,
\begin{align}
\pp [ N_{\Ieps} \geq k ] \leq N^{\delta} \eps^{ k (k+1) /2 }  
\end{align}

\eet

\proof  Our starting point is the inequality
\begin{align}
\pp [ N_{\Ieps} \geq 1 ] &\leq \eps^2 \ee [ \left( \Im [ m_N ( E + \i \eta )]\right)^2 ]\leq \frac{\eps^2}{N^2} \sum_{i, j} \ee [| G_{ii} | | G_{jj} | ] \leq \frac{\eps^2}{N^2} \sum_{i, j} \ee [ | G_{ii} |^2 ]^{1/2} \ee [  | G_{jj} |^2 ]^{1/2} \label{eqn:lr1}
\end{align}
where $N \eta = \eps$.  We will compute the expectations appearing on the right.

We have the formula
\beq
G_{ii} = \left( V_i + \sqrt{T} w_{ii} - E - \i \eta - \sum_j d_j \xi_j - \i \sum_j c_j \xi_j  \right)^{-1}
\eeq
where
\beq
c_j = \frac{ \eps} { N^2 ( \lami_j - E )^2 + \eps^2 } , \quad d_j = \frac{ N ( \lami_j - E ) }{ N^2 ( \lami_j - E )^2 + \eps^2 }
\eeq
and 
\beq
\xi_j =N \left( \sum_{k , k \neq i } h_{ik} u_{k}^{(j)} \right)^2
\eeq
where $h_{ik} = V_{i} \delta_{ik} + \sqT w_{ik}$ are the entries of $H_T$, $u^{(j)}$ is the $j$th normalized eigenvector of the $i$th minor of $H_T$ and $\lambda_j^{(i)}$ is the $j$th eigenvalue of the $i$th minor.   

From the formula
\beq
\frac{i}{a} = \int_0^\infty \e^{ \i u a } \d u , \qquad \Im [a] > 0 \label{eqn:lrb1}
\eeq
we have
\beq
| G_{ii} |^2 = \int_0^\infty \int_0^\infty \exp \left[ - ( u + v ) \sum_j c_j \xi_j  + \i ( u - v ) \sum_j d_j \xi_j + \i ( u  - v ) ( E - V_i - \sqT w_{ii} ) - ( u + v ) \eta  \right] \d u \d v.
\eeq
We denote by $\ee_i$ the expectation over the $i$th row of $W$. This is a Gaussian integral which we can compute explicitly, as conditionally on $H^{(i)}_T$, the variables $x_j = \sum_{k\neq i } N^{1/2} h_{i k } u_k^{(j)}$ are independent Gaussian random variables with variance $T$.  Hence,
\begin{align}
\ee_i [ | G_{ii} |^2 ] &= \int_0^\infty \int_0^\infty \int_{\rr^N}   \exp \bigg[  - \sum_{j=1}^{N-1} \frac{x_j^2}{2} (1 + 2 T(u + v)  c_j + \i (u - v) 2 T d_j  ) \notag \\
&- y^2 + \i ( u - v) ( E - V_i - \sqrt{2T} N^{-1/2} y ) - ( u + v ) \eta \bigg]  \frac{ \d y \d \boldsymbol{x} }{ ( 2 \pi )^{N/2}}\d u \d v \notag \\
&= \int_0^\infty \int_0^\infty \frac{ \e^{ - ( u - v )^2 \frac{ T}{2N} - ( u + v ) \eta + \i ( u - v ) (E - V_i ) }}{ \prod_l ( 1 + 2 T( u + v ) c_l + \i 2 T ( u - v ) d_l )^{1/2} } \d u \d v.
\end{align}
Note that we also integrated out the variable $w_{ii}$.  After the change of variables $2T(u +v, u-v) \to (u, v)$ this becomes
\beq
\ee_i [ |G_{ii} |^2 | ] = \frac{1}{ 4 T^2} \int_0^\infty \d u \int_{-\infty}^\infty \d v \frac{ \e^{ - v^2 \frac{ 1}{8 NT} - u \eta (2T)^{-1} + \i v ( E - V_i )  (2T)^{-1} } } { \prod_l ( 1 + u c_l + \i v d_l )^{1/2} } \label{eqn:lr2}.
\eeq

Let $1 > q_1 > q $. We let $\Gddi$ be the event
\begin{align} \label{eqn:lrf3}
\Gddi := \left\{ \left| \{ l : \lambda_l ( H^{(i)}_T ) \in \I_\eps \} \right| \leq N^{\delta/10} \right\} 
\end{align}
We know that $\Gddi$ holds with $(\xi, \nu)$-high probability.  Furthermore, $\Gddi$ is independent of the $i$th row and column of $H_T$.  Let us first get a bound for $\ee [ \onegddc |G_{ii}|^2]$.  On this event we know that there are at least $4$ eigenvalues in $\I_\eps$ and so there are at least $4$ indices $\alpha_j$ s.t. $c_{\alpha_j} \geq \eps^{-1}$.  Choosing such indices we bound
\begin{align}
\onegddc \ee_i [ |G_{ii}|^2 ] &\leq \frac{ C \onegddc }{T^2} \int_{0}^\infty \frac{1}{ \prod_{j=1}^3 | 1 + c_{\alpha_j} u |^{1/2} } \d u \times \int_{-\infty}^\infty \e^{ - v^2 / (NT) } \d v \notag \\
&\leq N^C \onegddc \eps
\end{align}
by our choice of $c_{\alpha_j}$.  Hence,
\beq
\ee [ \onegddc |G_{ii} |^2 ] \leq \eps N^C \e^{ - \nu \log (N)^{\xi}} \leq \eps. \label{eqn:lrf1}
\eeq

We now work on the event $\Gddi$.  Our goal is to obtain
\begin{align}
\ee [\onegdd  | G_{ii} |^2 ] \leq \frac{1}{ ( V_i - E)^2 + T^2 } N^{ \delta} \eps^{-1-r} . \label{eqn:fix1}
\end{align}
This will be the result of obtaining two bounds on the integral on the RHS of (\ref{eqn:lr2}) and taking the minimum at the end.  We begin the derivation of the first bound.  We choose distinct indices $\alpha_i$ and $\beta_i$ and use the bound
\begin{align} 
\onegdd \ee_i [ |G_{ii}|^2 ] &\leq \onegdd \frac{1}{4T^2} \int_0^\infty \d u \frac{1}{ | 1 + c_{\alpha_1} u |^{1/2} | 1 + c_{\alpha_2} u |^{1/2} | 1 + c_{\alpha_3} u |^{r}  } \notag \\
&\times  \int_{- \infty}^\infty \d v  \frac{1}{ ( 1 + |d_{\beta_1 } v |)^{1/2} ( 1 + |d_{\beta_2 } v |)^{1/2}( 1 + |d_{\beta_3 } v |)^{r}} . \label{eqn:lra3}
\end{align}
We now define a random variable $\Delone$ as follows.  On the complement of $\Gddi$ it is $0$.  If $\Gddi$ holds and there are at least $20$ eigenvalues of $H_T^{(i)}$ greater than $E + \eps/N$ we set it equal to $N | \lambda_l ( H_T^{(i)} ) - E |$ where $l$ is the index of the $20$th such eigenvalue (i.e., the smallest $l$ s.t. $\lambda_{l-20} ( H_T^{(i)} )$ is still greater than $E + \eps/N$).  If $\Gddi$ holds and there are less than $20$ eigenvalues greater than $E + \eps/N$ then there are at least $20$ eigenvalues less than $E- \eps /N$ and we set it equal to $N | \lambda_l ( H_T^{(i)} ) - E |$ where $l$ is the index of the $20$th such eigenvalue.  If $\Gddi$ holds and $\Delone <1$ we instead redefine it so $\Delone =1$ (this is just so $M \to \Delone^M$ is increasing).  By rigidity, the decay of the entries of $W$ and our assumption that $|| V|| \leq N^{B_V}$ it is easy to see that
\begin{align} \label{eqn:lrh2}
\ee [ \Delone^M ] \leq N^{\delta}
\end{align}
for any $M\geq 1$ and $\delta >0$.   By definition of $\Gddi$ we can always find distinct indices so that $c_{\alpha_i} \geq \eps/ (2 \Delone^2 )$ and $|d_{\beta_i}| \geq(2 \Delone)^{-1} $ (for example, take $\alpha_i$ so that $\lambda_{\alpha_i}$ are the closest eigenvalues to $E$ and then the $\beta_i$'s to be the next eigenvalues that are closest to $E$ but outside the interval $\Ieps$).  With this choice we get from (\ref{eqn:lra3}) that
\beq
\ee [ \onegdd  | G_{ii} |^2 ] \leq \frac{1}{ T^2} \eps^{-1-r} \ee [ \Delone^4] \leq\frac{1}{ T^2} \eps^{-1-r} N^{\delta}  .\label{eqn:lra2}
\eeq
For the other bound we start with the expression on the RHS of (\ref{eqn:lr2}) and integrate by parts in $v$ twice using
\beq
 \frac{-2T \i}{ E - V_i}\frac{ \d} { \d v} \e^{ \i v (E - V_i ) (2T)^{-1} } =  \e^{ \i v (E - V_i ) (2T)^{-1} }
\eeq 
both times.  Using $(N T )^{-1} \leq 1$, we bound the resulting integral by
\begin{align}
\onegdd \ee_i [ | G_{ii} |^2 ] &\leq \frac{C \onegdd}{ ( V_i - E )^2} \int_0^\infty \d u \int_{- \infty}^\infty \d v \frac{1}{ \prod_l |1 + u c_l + \i v d_l |^{1/2} } \notag  \\
&\times \left( 1 + | v|^2 + |v| \sum_j \frac{ | d_j | }{ |1 + u c_j + \i v d_j | } + \sum_{j, k} \frac{ | d_j | }{ |1 + u c_j + \i v d_j | } \frac{ | d_k | }{ |1 + u c_k + \i v d_k | }  \right).
\end{align}
Define the random set of indices
\beq
\F_{\delta} := \{ l : |\lambda_l ( H_T^{(i)} )- E | \leq N^{ \delta -1} \}. \label{eqn:fdelt}
\eeq
Note that by rigidity we have that
\begin{align}
\ee [ | F_{\delta} |^M ] \leq N^{2 \delta M} \label{eqn:lrh1}
\end{align}
for any $M \geq 1$.

By the Schwarz inequality we have
\begin{align}
& \onegdd \ee_i [ | G_{ii} |^2 ] \leq  \frac{C \onegdd }{ ( V_i - E )^2} \int_0^\infty \d u \int_{- \infty}^\infty \d v \frac{1}{ \prod_l |1 + u c_l + \i v d_l |^{1/2} }  \notag \\
&\times \left( 1 + |v|^2 + \sum_{ j, k \in \F_{\delta} } \frac{ | d_j | }{ |1 + u c_j + \i v d_j | } \frac{ | d_k | }{ |1 + u c_k + \i v d_k | } +  \sum_{ j, k \notin \F_{\delta} } \frac{ | d_j | }{ |1 + u c_j + \i v d_j | } \frac{ | d_k | }{ |1 + u c_k + \i v d_k | } \right) \notag \\
&=: A_1 + A_2 + A_3 + A_4 \label{eqn:lr3}
\end{align}
The argument above leading to the bound (\ref{eqn:lra2}) yields
\beq
\ee [ A_1 ] \leq \frac{1}{ ( V_i - E)^2} N^{  \delta } \eps^{-1-r}  \label{eqn:lrd1}
\eeq
i.e., we bound $A_1$ by the RHS of (\ref{eqn:lra3}) except that we have a prefactor of $(V_i - E)^{-2}$ instead of $T^{-2}$.  The term $A_2$ can be bounded similarly; to be more precise, it is bounded by the same thing as (\ref{eqn:lra3}) (with of course $T^{-2}$ replaced by $(V_i - E)^{-2}$) except we take  $4$ additional $d_{\beta_i}$'s for the denominator to deal with the $v^2$ in the numerator.  We get
\beq
\ee [ A_2 ] \leq \frac{1}{ ( V_i - E)^2} N^{ \delta } \eps^{-1-r}  \label{eqn:lrd2}
\eeq
We next handle $A_4$.  Define the random variable
\begin{align}
Q_\delta := \frac{1}{N} \sum_{j \notin \F_\delta} \frac{1}{ | \lami_j - E | }  \label{eqn:lrh3}.
\end{align}
We have the following bound which is a consequence of rigidity and the proof of Lemma \ref{lem:gilogbd}:
\beq
\ee [ Q_\delta^M ] \leq N^{2M \delta }\label{eqn:lr5}
\eeq
We provide the proof of this in Lemma \ref{lem:dbound} below.  Repeating the argument as before we then get that
\beq
\ee [ A_4 ] \leq  \frac{1}{ ( V_i - E )^2 } \ee [ Q_\delta^2 \Delone^4 ] \leq N^{10 \delta}
\eeq
where we have used the Schwarz inequality in the last step together with (\ref{eqn:lrh2}) and (\ref{eqn:lr5}).

The double sum in $A_3$ contains only $| F_\delta|^2$ terms.  Since we have (\ref{eqn:lrh1}) we will just bound each term individually.  We have,
\begin{align}
\onegdd \int_0^\infty \d u \int_{- \infty}^\infty & \d v \frac{1}{ \prod_l |1 + u c_l + \i v d_l |^{1/2} }  \frac{ | d_j | }{ |1 + u c_j + \i v d_j | } \frac{ | d_k | }{ |1 + u c_k + \i v d_k | } 
\notag\\
 \leq & \onegdd \int_0^\infty \d u \frac{1}{ ( 1 + c_{\alpha_1 } u )^{1/2} ( 1 + c_{\alpha_2} u )^{r} ( 1 + u c_j )^{1/4} ( 1 + u c_k )^{1/4} } \notag \\
&\times  \int_{-\infty}^\infty  \d v  \frac{ | d_j d_k | }{ ( 1 + |d_j v | )^{1/2} ( 1 + | d_k v | )^{1/2} ( 1 + | d_\beta v | )^r } =:A_{jk} \label{eqn:lr7}
\end{align}
Observe that due to the extra terms appearing in the denominator of the $d_j d_k$ terms of (\ref{eqn:lr3}) it is possible to make the choice $\alpha_1 = j$ or $k$ in the above inequality. This will be important later when we deal with the event $\{ N_{\I_\eps} \geq 2 \}$. Let us first do the $u$ integration.  The region of integration is split into two parts: $[0, \infty ) = [0, \eps^{20} ] \cup [ \eps^{20} , \infty )$. In the first region we bound the integrand by $1$.  The latter integral we bound by
\beq
\int_{\eps^{20}}^\infty \frac{1}{ c_{\alpha_1}^{1/2} c_{\alpha_2}^r c_j^{1/4} c_k^{1/4} u^{1 +r } } \d u \leq \frac{ \eps^{-20r}}{c_{\alpha_1}^{1/2} c_{\alpha_2}^r c_j^{1/4} c_k^{1/4} }.  \label{eqn:lrb5}
\eeq
Performing a similar split for the $v$ integration we get the bound
\beq
A_{jk} \leq ( \eps^{10} + \frac{ \eps^{-20r}}{c_{\alpha_1}^{1/2} c_{\alpha_2}^r c_j^{1/4} c_k^{1/4} }  ) ( \eps^{10} + \frac{ \eps^{ -20 r} |d_j|^{1/2} | d_k |^{1/2}  }{ | d_\beta|^r } ) \leq \eps^{5} N^{ \delta} \Delone^4 +  \eps^{-50r -1/2} \Delone^4\frac{ |d_j |^{1/2}}{  c_j^{1/4}  } \frac{ | d_k |^{1/2} } { c_k^{1/4}}. \label{eqn:lrb6}
\eeq
We have used the definition of $\Gddi$ to choose $c_{\alpha_i} \geq \eps / (2 \Delone^2 )$ and $| d_\beta| \geq(2  \Delone)^{-1}$.  Note also by the definition of $\F_\delta$ that $c_j, c_k \geq \eps N^{-\delta}$.  We have used this bound for $c_j$ and $c_k$, as well as the bound $|d_l| \leq \eps^{-1}$ which holds for any $l$, to bound the other terms that arise when the LHS of (\ref{eqn:lrb6}) is expanded (i.e., the terms of the form $\eps^{10} |d_j|^{1/2} |d_k|^{1/2} / |d_\beta|^r$, etc.).   Observe that
\beq
\frac{ |d_l |^2}{ c_l } \leq \frac{1}{\eps}
\eeq
for any $l$.  We have therefore derived the bound
\beq
\ee [ A_4 ] \leq N^{2 \delta } \eps^{-1-50r} \frac{1}{ ( V_i - E)^2 } \ee [ | F_\delta |^2 ( \Delone^4 ) ] \leq N^{8 \delta} \eps^{-1-50r} \frac{1}{ ( V_i - E)^2 } 
\eeq
where we have used the Schwarz inequality and (\ref{eqn:lrh2}), (\ref{eqn:lrh1}).

Collecting everything we get
\begin{align} \label{eqn:lrg1}
\ee [ \onegdd | G_{ii} |^2 ] \leq \frac{1}{ ( V_i - E)^2 + T^2 } N^{ 11 \delta} \eps^{-1-50r},
\end{align}
i.e., we have derived (\ref{eqn:fix1}).  As a consequence of (the proof of) Lemma \ref{lem:gilogbd} we have
\beq \label{eqn:lrg2}
\frac{1}{N} \sum_i \frac{1}{ \sqrt{ ( V_i - E)^2 + T^2 } } \leq N^{\delta}.
\eeq
Therefore after redefining $\delta$ and $r$ and using (\ref{eqn:lrf1}), (\ref{eqn:lrg1}) and (\ref{eqn:lrg2}) we get
\begin{align}
\frac{\eps^2}{N^2} \sum_{i, j} \ee [ | G_{ii} |^2]^{1/2} \ee [ |G_{jj} |^2 ]^{1/2} &\leq  C \eps^2\left( \eps^{1/2} + \frac{1}{N} \sum_{i=1}^N \frac{N^{\delta} \eps^{-1/2-r}}{ \sqrt{ (V_i - E)^2 + T^2} }  \right)^{2} \notag \\
&\leq C N^{ 4 \delta} \eps^{1-2r}
\end{align}
and so we obtain the bound in the case $k=1$.

We now begin the $k=2$ case.  The starting point is the inequality
\begin{align}
\pp [  N_{\Ieps} \geq 2 ] \leq \eps^2 \ee [  \1_{ \{ N_{\Ieps} \geq 2 \} } \left( \Im [ m_N ] \right)^2 ] \leq \frac{\eps^2}{N^2} \sum_{i, j} \ee [  \1_{ \{\Ni_{\Ieps} \geq 1 \}} | G_{ii} |^2 ]^{1/2} \ee [  \1_{ \{ \Nj_{\Ieps} \geq 1 \}} | G_{jj} |^2 ]^{1/2} \label{eqn:lrb7}
\end{align}
where $\Ni_{\Ieps}$ denotes the number of eigenvalues of the $i$th minor of $H_T$ lying in the interval $\Ieps$.  We have used the interlacing property which implies that if two eigenvalues lie in $\Ieps$ then at least one eigenvalue of any minor lies in $\Ieps$.  

Define the event $\Gddi$ as before.  Due to (\ref{eqn:lrf1}) we have
\beq
\ee [ \onegddc \oneni | G_{ii} |^2 ] \leq \eps. \label{eqn:lrf2}
\eeq
On the event $\Gddi$ we will prove the inequality
\beq
\ee [  \onegdd \1_{ \{ \Ni_{\Ieps} \geq 1 \} }  |G_{ii} |^2 ] \leq N^{\delta} ( (V_i - E)^2 + T^2 )^{-1} \eps^{-r} \left( \pp [ \Ni_{\Ieps } \geq 1 ]\right)^{1/(1+r)}. \label{eqn:lr8}
\eeq
Proceeding in the $k=1$ case we use the identity (\ref{eqn:lrb1}) and take the expectation over the $i$th row.  Note that $\Ni_{\Ieps}$ is independent of the $i$th row of $H_T$.  We see that
\beq
\onegdd \ee_i [ | G_{ii} |^2 \1_{ \{ \Ni_{\Ieps } \geq 1 \} } ] =  \onegdd \oneni \frac{1}{ 4 T^2} \int_0^\infty \d u \int_{-\infty}^\infty \frac{ \e^{ - v^2 \frac{ 1}{ 8 N T} - u \eta (2T)^{-1} + \i v ( E - v_i ) (2T)^{-1} } } { \prod_l ( 1 + u c_l + \i v d_l )^{1/2} }. \label{eqn:lrb2}
\eeq
Again, we get two bounds on the integral and minimize over the two at the end. The arguments are very similar to the $k=1$ case, except that \emph{due to the fact that $\Ni_{\Ieps} \geq 1$, we will always be able to choose an index $\alpha_1$ so that $c_{\alpha_1} \geq \eps^{-1}$}.    In the case $k=1$ we could only choose $c_{\alpha_1} \geq \eps/ (2 \Delone^2 ) \sim  \eps N^{- \delta}$ by rigidity.  This was one of the key observations in the work \cite{wegner}.

For the first inequality, we proceed as in (\ref{eqn:lra3}) and derive
\begin{align}
\onegdd \oneni \ee_i [ |G_{ii}|^2 ] &\leq  \onegdd \oneni \frac{C}{T^2} \int_0^\infty \d u \frac{1}{ | 1 + c_{\alpha_1} u |^{1/2} | 1 + c_{\alpha_2} u |^{1/2} | 1 + c_{\alpha_3} u |^{r}  } \notag \\
&\times 
\ \int_{- \infty}^\infty \d v  \frac{1}{ ( 1 + |d_{\beta_1 } v |)^{1/2} ( 1 + |d_{\beta_2 } v |)^{1/2}( 1 + |d_{\beta_3 } v |)^{r}} .
\end{align}
Since at least one eigenvalue of the minor lies in $\Ieps$ so we can choose $c_{\alpha_1} \geq \eps^{-1}$.  As before, we then choose the other indices so that $c_{\alpha_2} \geq \eps / ( 2 \Delone^2 )$ and $|d_{\beta_i} | \geq (2 \Delone)^{-1}$ which is possible by the definition of $\Gddi$. Hence,
\beq \label{eqn:lrb3}
\ee [  \onegdd \oneni | G_{ii} |^2 ] \leq \frac{1}{T^2}  \eps^{-r} \ee [ \Delone^4 \oneni ] \leq \frac{1}{T^2} N^{8 \delta} \eps^{-r} \pp[ \Ni_{\Ieps} \geq 1 ]^{1/(1 +r)}
\eeq
where we have used the H\"older inequality in the last step and (\ref{eqn:lrh2}).  
 For the other bound we proceed as before and integrate (\ref{eqn:lrb2}) by parts twice in $v$.  We obtain
\beq
\onegdd \oneni  \ee_{i} [ | G_{ii} |^2 ] \leq \oneni (A_1 + A_2 + A_3 + A_4 ) =: B_1 + B_2 + B_3 + B_4
\eeq
where the $A_i$ are defined as in (\ref{eqn:lr3}).  The same argument leading to (\ref{eqn:lrb3}) (i.e., the same arguments we used to obtain (\ref{eqn:lrd1}) and (\ref{eqn:lrd2}) except taking $c_{\alpha_1} \geq \eps^{-1}$)  yields
\beq
\ee [ B_1 + B_2 ] \leq \frac{C}{(V_i - E)^2} \eps^{-r} \ee [ \Delone^6 \oneni ] \leq  \frac{C}{(V_i - E)^2} \eps^{-r} N^{12 \delta} \pp[\Ni_{\Ieps} \geq 1 ]^{1/(1+r)}.
\eeq
Similarly we get
\begin{align}
\ee [ B_4 ] &\leq \frac{1}{ ( V_i - E )^2} \eps^{-r} \ee [ Q_\delta^2 \Delone^4 \oneni ] \notag \\
&\leq \frac{1}{ ( V_i - E )^2} N^{12 \delta} \eps^{-r} \pp [ \Ni_{\Ieps } \geq 1]^{1 / (1+r)}  
\end{align}
using (\ref{eqn:lr5}), (\ref{eqn:lrh2}) and the H\"older inequality. We proceed as before in our estimation of $B_3$.  We estimate the general term by 
\begin{align}
 \onegdd \oneni \int_0^\infty \d u \int_{- \infty}^\infty & \d v \frac{1}{ \prod_l |1 + u c_l + \i v d_l |^{1/2} }  \frac{ | d_j | }{ |1 + u c_j + \i v d_j | } \frac{ | d_k | }{ |1 + u c_k + \i v d_k | } 
 \notag \\
 \leq &  \onegdd \oneni \int_0^\infty \d u \frac{1}{ ( 1 + c_{\alpha_1 } u )^{1/2} ( 1 + c_{\alpha_2} u )^{r} ( 1 + u c_j )^{1/4} ( 1 + u c_k )^{1/4} } \notag \\
&\times  \int_{-\infty}^\infty  \d v  \frac{ | d_j d_k | }{ ( 1 + |d_j v | )^{1/2} ( 1 + | d_k v | )^{1/2} ( 1 + | d_\beta v | )^r } =:B_{jk}  \label{eqn:lrb4}
\end{align}
Note that we can take $\alpha_i = j$ or $k$ due to the extra factors appearing in the denominator below the $d_j$ and $d_k$ terms in the first line of (\ref{eqn:lrb4}) (this was the observation made after (\ref{eqn:lr7})).  Due to the fact that $\Ni_{\Ieps} \geq 1$ we can take $c_{\alpha_1} \geq \eps^{-1}$.  By the definition of $\Gddi$ we can take $c_{\alpha_2} \geq \eps/ (2 \Delone^2)$ and $b_{\alpha} \geq (2\Delone)^{-1}$.  Splitting the two integrals into the regions $[0, \eps^{20} ] \cup [\eps^{20} , \infty )$ as in (\ref{eqn:lrb5}) we get, in the same fashion as (\ref{eqn:lrb6}), the bound
\begin{align}
B_{jk} &\leq  \onegdd \oneni ( \eps^{10} + \frac{ \eps^{-20r}}{c_{\alpha_1}^{1/2} c_{\alpha_2}^r c_j^{1/4} c_k^{1/4} }  ) ( \eps^{10} + \frac{ \eps^{ -20 r} |d_j|^{1/2} | d_k |^{1/2}  }{ | d_\beta|^r } ) \notag \\
&\leq  \oneni \left( \eps^{5} N^{4 \delta} \Delone^4 +  \eps^{-50r +1/2} \Delone^4\frac{ |d_j |^{1/2}}{  c_j^{1/4}  } \frac{ | d_k |^{1/2} } { c_k^{1/4}} \right) \notag \\
&\leq \eps^{-50r} \Delone^4 \oneni.
\end{align}
From H\"older and (\ref{eqn:lrh2}), (\ref{eqn:lrh1}) we get
\begin{align}
\ee [ B_3 ] \leq \eps^{-50r} N^{4 \delta} \ee [ \Delone^4 | F_\delta|^2 \oneni ] \leq \eps^{-50r} N^{ 12 \delta} \pp [ \Ni_{\Ieps} \geq 1]^{1/(1+r)}.
\end{align}
Collecting everything and optimizing over our two bounds we get
\begin{align}
\ee [ \onegdd \oneni | G_{ii} |^2 ] &\leq \frac{1}{ ( V_i - E)^2 + T^2} N^{ \delta} \eps^{-r} \pp [ \Ni_{\Ieps} \geq 1 ]^{1/(1+r)}  \notag \\
&\leq  \frac{1}{ ( V_i - E)^2 + T^2} N^{\delta} \eps^{1-3r}  
\end{align}
where in the last step we have applied the $k=1$ result to the $i$th minor of $H_T$.  Plugging this and (\ref{eqn:lrf2}) into (\ref{eqn:lrb7}) and summing over $i, j$ yields the claim.  \qed

Outside the interval $\I_{E_0, qG}$ we still have the following bound.  It is a corollary of the above proof.  It may also be established using the methods of \cite{homog, wegner}.

\bec \label{cor:lrest}
For all $r>0$ there is a $C_r >0$ so that
\begin{align}
\pp [ N_{\Ieps } \geq k ] \leq N^{C_r} \eps^{k(k+1)/2} ( |E - E_0 | + 1)^{C_r}
\end{align}
for $k=1, 2$.
\eec
\proof The proof is essentially the same as the proof of Theorem \ref{thm:lrest}.  We illustrate the differences in the case $k=1$.   Since now the event $\Gddi$ may have probability $1$, the bound (\ref{eqn:lrf1}) becomes
\begin{align}
\ee [ \onegddc |G_{ii} |^2 ] \leq \eps N^C.
\end{align}
On the event $\Gddi$ there is no need to integrate by parts as we do not need to optimize the prefactor.  We just use the argument leading to (\ref{eqn:lra2})  and get
\begin{align}
\ee [ \onegdd  | G_{ii} |^2 ] \leq \frac{1}{ T^2} \eps^{-1-r} \ee [ \Delone^4] \leq N^2 \eps^{-1-r} \ee [ \Delone^4]
\end{align}
Since all the eigenvalues of $H_T^{(i)}$ lie in the interval $[ A_V, N^{B_V}] + [ - ||W^{(i)} ||, ||W^{(i)} ||]$ it is easy to see that
\begin{align}
\ee [ \Delone^M ] \leq N^{CM} ( | E - E_0| +1 )^{CM}.
\end{align}
This gives the $k=1$ bound.  The $k=2$ bound is similar. \qed

\noindent{\bf Proof of Theorem \ref{thm:lr}}.  We decompose the event $\{ | \lamit - \lamipt | \leq \eps / N \}$ as follows.
\begin{align}
\{ | \lamit - \lamipt | \leq \eps / N \} &= \{  | \lamit - \lamipt | \leq \eps / N , | \lamit - \gamit | \leq N^{\delta} / N \} \notag \\
& \bigcup_{n=1}^\infty \{  | \lamit - \lamipt | \leq \eps / N , | \lamit - \gamit | \geq N^{\delta} / N , n-1 \leq || H_T || \leq n \} \\
& =: \F_1 \cup \bigcup_n \G_n.
\end{align}
Define
\beq
\I_j := [ \gamit + \eps (j-2) / N , \gamit + \eps (j+2) / N ].
\eeq
Using Theorem \ref{thm:lrest} we have
\beq
\pp [ \F_1 ] \leq \sum_{ |j| \leq N^{ \delta} / \eps } \pp [ N_{\I_j } \geq 2 ] \leq N^{3 \delta} \eps^{2 -r }.
\eeq
Let $p_1, p_2$ and $p_3$ satisfy $p_1^{-1} + p_2^{-1} + p_3^{-1} = 1$ with $p_1 > 1$. We have by H\"older's inequality and Corollary \ref{cor:lrest}
\begin{align}
\pp [ \G_n ] & \leq \sum_{|j| \leq 4n / \eps } \pp [ N_{\I_j} \geq 2 , | \lamit - \gamit | \geq N^{\delta} / N , n-1 \leq || W || \leq n ] \notag \\
&\leq N^{C_r} \eps^{ (3-r)/p_1 -1} n^{C_r} \left(\pp [ ||W || \geq n-1] \right)^{1/p_2} \e^{ - \nu \log (N)^{\xi} / p_3 }.
\end{align}
For large enough $N$ depending on $p_2$ and $p_3$ and $C_r$ we have due to the exponential decay of the entries of $W$,
\beq
\sum_{n=1}^\infty N^{C_r} n^{C_r} \left( \pp [ ||W || \geq n-1 ] \right)^{1/p_2} \e^{ - \nu \log (N)^{\xi}/ p_3}  \leq 1.
\eeq
Hence choosing $p_1$ close to $1$ we obtain the theorem. \qed

We required the following lemma in the proof of Theorem \ref{thm:lrest}.
\bel \label{lem:dbound} We have,
\begin{align}
\ee [ Q_\delta^M] \leq N^{ 2M \delta}.
\end{align}
for $E \in \I_{E_0, q G}$ and $\F_{\delta}$ as in (\ref{eqn:fdelt}) and $Q_\delta$ as in (\ref{eqn:lrh3}).
\eel
\proof Let $1 > q_1 > q > 0$ and define $\Qone$ to be the event that $| \gamma_{j, T } - \lami_j | \leq N^{\delta/10}$ for $j \in \A_{q_1, T}$ and $|| W^{(i)} || \leq 3$.  Then $\Qone$ holds with $(\xi, \nu)$-high probability.  On the complement of $\Qone$ we just use the bound $|Q_\delta| \leq N^{2M}$ and so we may assume that $\Qone$ holds.  We write
\beq
\1_{\Qone} \sum_{j \notin \F_{\delta}} \frac{1}{N | \lami_j - E | } = \1_{\Qone} \sum_{\substack{ j \notin \F_{\delta} \\ \gamjl \in \I_{E_0, q_1 G}}} \frac{1}{N | \lami_j - E | }  +  \1_{\Qone} \sum_{\substack{ j \notin \F_{\delta} \\ \gamjl \notin \I_{E_0, q_1 G}}} \frac{1}{N | \lami_j - E | } =: A_2 + A_1.
\eeq
For the terms in $A_2$ we get directly from the definitions of $\Qone$ and $\F_{\delta}$ that
\beq
A_2 \leq C \sum_{|k| \leq N} \frac{1}{|k|} \leq N^{2\delta}
\eeq
Fix a $q_2$ satisfying $q_1 > q_2 > q$.   Define the intervals for $k \geq 0$
\beq
\Ione_k = E_0 + [ -2^k G q_2, 2^k G q_2 ]
\eeq
as in the proof of Lemma \ref{lem:gilogbd}.  By definition, on the event $\Qone$ the norm of $W^{(i)}$ is less than $3$.  We have the bound $| \lambda_j (A) - \lambda_j (B) | \leq || A - B||$ which holds for matrices $A$ and $B$.  Hence we see that on the event $\Gddi$, 
\beq
| \lambda_j (H^{(i)}_T) - \lambda_j (V ) | \leq \sqT.
\eeq
Therefore, for $k \geq 1$ we see that
\beq
| \{ i : \lambda_i (H^{(i)}_T) \in \Ione_{k} \backslash \Ione_{k-1} \} | \leq  | \{ i : V_i \in \Ione_{k+1} \} | \leq C N 2^k G
\eeq
where the last inequality is proven in Lemma \ref{lem:gilogbd} and holds as long as $k \leq \left| \log (cG) / \log (2) \right|$ for some $c>0$.  From this we see that
\beq
A_1 \leq \frac{1}{N} \sum_{k=1}^{ C \log (G) } \frac{ C N 2^k G}{ 2^k G} \leq C \log (N) \leq N^{\delta}.
\eeq
The claim follows.
\qed

\section{Proofs of main results} \label{sect:mr}

\noindent{\bf Proof of Theorem \ref{thm:gap}. }  Fix $t= N^{\tau} / N   \in \Tom$ and $k_0 \in \Amut$. Let 
\beq
t_0 = t - \frac{N^{\om'}}{N}
\eeq
for $\om' < \tau / 100$.  We apply Theorem \ref{thm:matching} with the matching GOE matrix $\hatH_{t_0} = a \tilW + b$ for $a$, $b$ and $j_0$ chosen as in Lemma \ref{lem:matching}.
It follows that there is an $\eps >0$ and an event $\F$ satisfying $\pp [ \F ] \geq 1 - N^{-\eps}$ on which 
\beq \label{eqn:gap1}
\left|  \hatlam_{j+k_0, t} - \hatlam_{j'+k_0, t} - ( \hatmu_{j+j_0, t} - \hatmu_{j'+j_0, t} ) \right| \leq \frac{N^{-\eps}}{N},
\eeq
for $|j  | + |j' | \leq N^{c}$ for some $c>0$.  We take $c < \om / 20$.  Recall that $\bhatlam$ are distributed as $\blam ( H_t)$ and $\bhatmu$ are distributed as $\blam (a_t \tilW + b_t )$.  Clearly $a_t$ satisfies $|a - a_t | \leq | t  - t_0 |$.  

By Lemma \ref{lem:timeder},
\begin{align}
| \rhofcto ( \gamkoto ) - \rhofct ( \gamma_{k_0, t} ) | &\leq |  \rhofcto ( \gamkoto ) -  \rhofcto ( \gamma_{k_0, t} ) ) | + |  \rhofcto ( \gamma_{k_0, t} ) -  \rhofct ( \gamma_{k_0, t}  ) | \notag\\
& \leq \log (N) \frac{|t_0 - t|}{t_0} \leq \frac{1}{N^{\om/2}}
\end{align}
Since the eigenvalue differences appearing in (\ref{eqn:gap1}) are bounded by $N^{2c}/N \ll N^{ \om/4}/N$ with very high probability, we see that there is an event $\F_2$ with $\pp [ \F_2 ] \geq 1- N^{-\eps /2}$ on which 
\beq
\left| \rho_{\mathrm{fc}, t} ( \gamma_{k_0, t} ) (  \hatlam_{j, t} - \hatlam_{j', t} ) - \frac{\rhosc ( \mu_{j_0} )}{a_t} (\hatmu_{j, t} - \hatmu_{j', t} ) \right| \leq \frac{N^{-\eps}}{N}
\eeq
for $| j - k_0 | + | j' - k_0 | \leq N^{c}$, after possibly decreasing $\eps$. Now note that $a_t^{-1} \bhatmu$ is distributed as the eigenvalues of a standard GOE matrix.  The claim follows by expanding the test function $O$ on the set $\F_2$. \qed

\noindent{\bf Proof of Theorem \ref{thm:bulk}.}  The above proof establishes the fact that there is a $c > 0$ so that for any index $i \in \Amut$ and positive integers $k_1, ..., k_{n} \leq N^{c}$ we have
\begin{align}
\big| \ee^{H_t} &[ O ( \rhofct ( \gamit ) ( \lambda_{i} - \lambda_{i + k_1 } ) , ..., \rhofct ( \gamit) ( \lambda_{i} - \lambda_{i + k_n } ) ) ]  \notag \\
&- \ee^{GOE} [ O ( \rhosc ( \mu_i ) ( \lambda_{i} - \lambda_{i + k_1 } ) , ..., \rhosc ( \mu_i) ( \lambda_{i} - \lambda_{i + k_n } ) ) ] \big| \leq N^{-c}.
\end{align}
The method of going from the above estimate and optimal rigidity results to averaged bulk universality is standard and can be found for example in Section 7 of \cite{localrelaxation}.  We omit the details and only mention that our result is stated with an upper bound on the size of the averaging window which ensures that
\beq
| \rhofct (E) - \rhofct (E') | \leq N^{- c}
\eeq
for all energies satisfying $|E' - E | \leq b$.  If a larger averaging window is desired then $\rhofct (E)$ must be replaced by $\rhofct (E')$ in the theorem statement.
 \qed


\section{Proof of local deformed semicircle law}
\label{sect:llproof}
\subsection{Stability estimates}

In this section we will derive various qualitative properties of the deformed semicircle law.  All statements and proofs here are deterministic.  It will be notationally convenient to work with the definition
\beq
H_T= V + \sqT W
\eeq
and
\beq
\mfcl (z) = \frac{1}{N} \sum_{i=1}^N \frac{1}{ V_i - z - T \mfcl (z) } \label{eqn:mfcl}
\eeq
where we have dropped the subscript $T$ in $m_{\mathrm{fc}, T}$.  Throughout this section we will assume that $V$ is $(\ell, G)$-regular at the point $E_0$.

We will need the following elementary lemma about Stieltjes transforms.
\bel \label{lem:tilm}
Let $\tilm$ be the Stieltjes transform of a probability measure $\mu$ with support in $[-A_*, A_* ]$.  There is a constant $C >0$ s.t.
\beq \label{eqn:imtilm}
| \tilm (E + \i \eta ) | \leq C  \left( | \log ( \eta ) | + | \log (A_* ) | \right) \sup_{ \eta' \geq \eta}  \Im [ \tilm (E+ \i \eta' ) ] .
\eeq
\eel
\proof Define the dyadic intervals for $j \geq 0$ by
\beq
\I_j := [ E - 2^j \eta , E + 2^j \eta ]
\eeq
and set $\I_{-1} = \emptyset$.  We have
\beq
\Im [ \tilm (E + 2^j \eta ) ] \geq \int_{I_j} \frac{ 2^j \eta }{ (x - E)^2 + (2^j \eta )^2 } \d \mu (x) \geq c \frac{ \mu ( \I_j) }{2^j \eta }.
\eeq
Let $j_* = \sup_j \{ j : \mu ( I_j \backslash I_{j-1} )  \neq 0 \}$.  Then $j_* \leq C (| \log ( \eta ) | + | \log (A_* ) | )$.  We have
\beq
| \tilm (E + \i \eta ) | \leq \sum_{j=0}^{j_*} \int_{ \I_j \backslash \I_{j-1}} \frac{ \d \mu (x) }{ |x - z | } \leq C \sum_{j=0}^{j_*} \int_{ \I_j \backslash \I_{j-1} } \frac{ \d \mu (x) }{ 2^j \eta } \leq C \sum_{j=0}^{j_*} \frac{ \mu ( \I_j ) } { 2^j \eta }  \leq C j_* \sup_{\eta' \geq \eta } | \tilm (E + \i \eta' ) |.
\eeq
where we applied \eqref{eqn:imtilm} in the last inequality. The claim follows. \qed

\bel \label{lem:rstab} Assume that $V$ is $(\ell, G)$-regular at $E_0$. For  $1> q >0$ and $\om >0$ all of the following hold for $E \in \I_{E_0, qG}$ and $10 \geq \eta \geq N^{-5}$ and $T \in \Tom$, for $N$ large enough.  Below the constants do not depend on the choice of  $\om$ or $q$.  For the Stieltjes transform,
\beq \label{eqn:stab1}
c \leq \Im [ \mfcl (z) ] \leq C.
\eeq
and therefore
\beq \label{eqn:stab2}
c T \leq | V_i - z - T \mfcl (z) |\leq C.
\eeq
For the density $\rhofcl$ we have
\beq \label{eqn:stab3}
c \leq \rhofcl (E) \leq C , \qquad | \rhofcl ' (E) | \leq \frac{C}{T}.
\eeq
Finally,
\beq \label{eqn:stab4}
c \leq | 1 - T R_2 (z) | \leq C , \qquad | T^2 R_3 (z) | \leq C
\eeq
where 
\beq
R_k (z) := \frac{1}{N} \sum_i \frac{1}{ (v_i - z - T \mfcl (z) )^k} .
\eeq

\eel
\proof
By the assumptions on $V$ and Lemma \ref{lem:tilm} we see that there is a constant $C_V >1$ s.t. 
\beq \label{eqn:mvabs}
| m_V ( E + \i \eta ) | \leq C_V \log (N), \qquad E \in \I_{E_0, G}, \mbox{ } \eta \geq T / (10 C_V )
\eeq
and 
\beq \label{eqn:mvim}
C_V^{-1} \leq \Im [ m_V (E + \i \eta ) ] \leq C_V, E \in \I_{E_0, G} \mbox{ } 10 \geq \eta \geq T / (10 C_V ).
\eeq
For notational simplicty we define the set
\beq
\A = \left\{ E + \i \eta : E \in \I_{E_0, G} , \mbox{ } 10 \geq \eta \geq T / (10 C_V ) \right\}.
\eeq
Fix $E \in \I_{E_0, q G}$ and define
\beq \label{eqn:etastardef}
\eta_* = \inf \left\{ \eta \leq 5 : | \mfc (E + \i \eta ) | \leq 2 C_V \log (N), \mbox{ and } (2 C_V)^{-1} \leq \Im [ \mfc (E + \i \eta )] \leq 2 C_V \right\}.
\eeq
Using the bound $|\mfc | \leq 1$ for $\eta \geq 1$ we see that $\eta_* \leq 1$.  For a contradiction suppose that $\eta_* > 0$.  Since $\mfc$ is continuous in the upper half plane and the set defining $\eta_*$ in \eqref{eqn:etastardef} is nonempty (it has at least all $\eta \geq $) one of the following three possibilities must hold.
\begin{enumerate}[label=(\roman*)]
\item $| \mfc (E + \i \eta_* ) | = 2 C_V \log (N)$ and $(2 C_V )^{-1} \leq \Im [ \mfc (E + \i \eta_* ) ] \leq 2 C_V$ \label{item:s1}
\item $ | \mfc (E + \i \eta_*)| \leq 2 C_V \log (N)$ and $(2 C_V )^{-1} = \Im [ \mfc (E + \i \eta_* ) ]$ \label{item:s2}
\item $ | \mfc (E + \i \eta_*) |\leq 2 C_V \log (N)$ and $2 C_V = \Im [ \mfc (E + \i \eta_* ) ]$. \label{item:s3}
\end{enumerate}
We have the equation
\beq \label{eqn:mfcstar}
\mfc (E + \i \eta_* ) = m_V ( E + \i \eta_* + T \mfc (E + \i \eta_* ) ).
\eeq
Suppose that \ref{item:s1} holds. For $N$ large enough depending only on $C_V$, we have that $E + \i \eta_* + T \mfc (E + \i \eta_* ) \in \A$.  But then by \eqref{eqn:mfcstar} and \eqref{eqn:mvabs} we have
\beq
| \mfc (E + \i \eta_* ) | \leq C_V \log (N)
\eeq
which is a contradiction.  If \ref{item:s2} holds then still $E + \i \eta_* + T \mfc (E + \i \eta_* ) \in \A$.  By \eqref{eqn:mfcstar} and \eqref{eqn:mvim} we get
\beq
\Im [ \mfc (E + \i \eta_* ) ] \geq (C_V)^{-1}
\eeq
which is a contradiction.  If \ref{item:s3} holds we arrive at a similar contradiction.  Therefore, $\eta_* = 0$.  We conclude (\ref{eqn:stab1}),  (\ref{eqn:stab2}) and the first bound of (\ref{eqn:stab3}).

 Taking imaginary parts on both sides of (\ref{eqn:mfcl}) and rearranging we obtain
\beq
\Im [ \mfcl (z) ] \left( 1 - \frac{T}{N} \sum_{i=1}^N \frac{1}{ | V_i - z - T \mfcl (z) |^2 } \right) = \eta \left( \frac{1}{N} \sum_{i=1}^N \frac{1}{ | V_i - z - T \mfcl (z) |^2} \right) \label{eqn:blah1}
\eeq
and so
\beq
1 -\frac{T}{N} \sum_{i=1}^N \frac{1}{ | V_i - z - T \mfcl |^2 }  \geq 0. \label{eqn:s1}
\eeq

Denote $z + T \mfcl (z) = a + b \i$ for $a, b \in \rr$. We have,
\beq \label{eqn:staba1}
\Re [ 1 - T R_2 (z) ] = 1 - \frac{T}{N}\sum_{i=1}^N \frac{ 1}{ | V_i - a - b \i |^2 } + 2 b^2 \frac{T}{N} \sum_{i=1}^N \frac{1}{ | V_i - a - b \i |^4}.
\eeq
By the hypothesis on $V$ and the fact that $ a \in \I_{E_0, \L}$ we see that there is a $C >0$ so that if $ b \geq C T$ then
\beq
1 -\frac{T}{N}\sum_{i=1}^N \frac{ 1}{ | V_i - a - b \i |^2 }  \geq \frac{1}{2}.
\eeq
 We have by Lemma \ref{lem:vbound} below that there are at least $NT$ entries of $V$ less than distance $CT$ from $a$.  Hence in the regime $ C T \geq b \geq c T$ we have
\beq
b^2 \frac{T}{N} \sum_{i=1}^N \frac{1}{ | V_i - a - b \i |^4} \geq  b^2 \frac{T}{N} \frac{NT}{ (CT)^4}  \geq c
\eeq
and we conclude the first part of (\ref{eqn:stab4}), using (\ref{eqn:s1}) (note that by (\ref{eqn:stab1}) we always $b \geq cT$).  The second part of (\ref{eqn:stab4}) is an obvious consequence of (\ref{eqn:vass}).  The last thing to prove is the second bound of (\ref{eqn:stab3}).  Differentiating (\ref{eqn:mfcl}) gives
\beq
\del_z \mfcl (z) = \left( 1 - \frac{T}{N} \sum_{i=1}^N \frac{1}{ ( V_i - z - T \mfcl (z) )^2 }  \right)^{-1} \frac{1}{N} \sum_{i=1}^N \frac{1}{ ( V_i - z - T \mfcl (z) )^2} .
\eeq
Therefore, by the first bound of (\ref{eqn:stab4}), 
\begin{align}
| \del_z \mfcl (z) | \leq C \frac{1}{N} \sum_{i=1}^N \frac{1}{ | V_i - z - T\mfcl (z) |^2 }  \leq \frac{C}{ T}
\end{align}
where in the second inequality we have used (\ref{eqn:vass}) (note that we have used the fact that $| \Re [ \mfcl ] | \leq T^{-1/2}$ which implies that if $E \in \I_{E_0, q G}$ then $E + T\Re [ \mfcl ( E + \i \eta ) ] \in \I_{E_0, \sqrt{q}G}$).  We conclude the bound for $\rhofc'$ by the Stieltjes inversion formula. \qed

\bel \label{lem:dumbstab}
For $z \in \D_2$ where $\D_2$ is defined in (\ref{eqn:donedtwo}) we have,
\beq \label{eqn:dumbstab}
 | V_i - z - T \mfc (z) | \geq c, \quad c \leq | 1 - T R_2 (z) | \leq C , \quad | T^2 R_3 (z) | \leq C
\eeq
and
\beq
|  \mfc (z)  |\leq C \label{eqn:dumbmbd}
\eeq

\eel
\proof  
The first bound follows from the fact that $\Im [z] \geq 10$ on $\D_2$.
The other two bounds are consequences of the fact that we have $|R_k| \leq C$ for $z \in \D_2$.  The bound (\ref{eqn:dumbmbd}) is obvious. \qed

For the above proofs we required the following elementary result.
\bel \label{lem:vbound} Fix $1 >q>0$ and $\om > 0$, and $V$ be $(\ell, G)$-regular at $E_0$.  Then for $N$ large enough we have uniformly for $E \in \I_{E_0 , q G}$ and $ \eta \in \Tom$,
\beq \label{eqn:vbound}
c N\eta  \leq \left| \{ i : V_i \in [E - \eta, E + \eta ] \} \right| \leq C N \eta
\eeq
for some constants depending only on the constants in (\ref{eqn:vass}).

\eel

\proof Let $\Ieta = [E - \eta, E + \eta ] $.  The upper bound follows immediately from the fact that
\beq
V_i \in \Ieta \implies ( V_i - E )^2 + \eta^2 \leq 2 \eta^2.
\eeq
Denote by $C^*$ the constant appearing in the upper bound just obtained for (\ref{eqn:vbound}).  For $k \in \zz$ denote
\beq
\Iketa = \eta k + \Ieta,
\eeq
and define
\beq
k_0 = \inf \left\{ |k| : \Iketa \cap \{ \rr \backslash \I_{E_0, \sqrt{q} G } \} \neq \emptyset \right\}.
\eeq
We have,
\begin{align}
\Im [ \mV (E + \i \eta )] \leq \sum_{|k| \leq k_0} \frac{1}{N} \sum_{V_i \in \Iketa} \frac{\eta}{ ( V_i - E)^2 + \eta^2 } + C \frac{\eta}{G^2}.
\end{align}
Note that $C \eta / G^2 \leq C N^{-\om}$ by assumption on $\eta$.
Fix $L > 0$, $L \in \nn$ and let $\NL := | \{ i : V_i \in [E - L\eta, E + L\eta ] \}|$.  We have,
\begin{align}
\sum_{|k| \leq k_0 } \frac{1}{N} \sum_{v_i \in \Iketa} \frac{\eta}{ ( V_i - E)^2 + \eta^2 }  & \leq \sum_{|k| \leq L} \frac{1}{N} \sum_{V_i \in \Iketa} \frac{\eta}{ ( V_i - E)^2 + \eta^2 }  + \sum_{k_0 \geq |k| > L} \frac{1}{N} \sum_{V_i \in \Iketa} \frac{\eta}{ ( V_i - E)^2 + \eta^2 } \notag  \\
& \leq \frac{\NL}{N \eta} + \sum_{ |k| > L} \frac{ 4 C^* }{ k^2} \leq \frac{ \NL}{ N \eta} + \frac{ 20 C^*}{ L}.
\end{align}
Choosing $L > 100 C^* / c_V$ where $c_V$ is the lower bound from (\ref{eqn:vass}) we get for large $N$,
\beq
\NL \geq ( N \eta ) \frac{c}{2}
\eeq
which yields the claim. \qed

Define
\beq
g_i (z) = \frac{1}{ V_i - z - T \mfcl (z) }.
\eeq
We will require the following bound later.

\bel \label{lem:gilogbd} For all $z \in \D_{L, q}$, we have
\beq
\frac{1}{N} \sum_i |g_i (z) | \leq C \log N.
\eeq
\eel
\proof This is clear for $z \in \D_2$.  For $z \in \D_1$ we write
\begin{align}
\frac{1}{N} \sum_i |g_i (z) | = \frac{1}{N} \sum_{i : V_i \in \I_{E_0, \sqrt{q}G} } |g_i (z) | + \frac{1}{N} \sum_{i : V_i \notin \I_{E_0, \sqrt{q} G} } |g_i (z) | =: A_2 + A_1.
\end{align}
We first bound $A_1$. Define for $k \geq 0$
\beq
\Ione_k := E_0 + [ -2^k \sqrt{q}G, 2^k \sqrt{q}G], \qquad N^{(1)}_k := | \{ i : V_i \in \Ione_k \} |.
\eeq
For $k \leq k_1 := \left\lfloor{-\log (G) / \log (2) + 1} \right\rfloor$ we have by the assumptions (\ref{eqn:vass}) and the proof of Lemma \ref{lem:vbound} the bound
\beq
N^{(1)}_k \leq C 2^k G N.
\eeq
Hence,
\begin{align}
A_1 &= \frac{1}{N} \sum_{1 \leq k \leq k_1 } \sum_{V_i \in \Ione_k \backslash \Ione_{k-1} } |g_i (z) | + \frac{1}{N} \sum_{V_I \notin \Ione_{k_1} } |g_i (z) | \notag \\
&\leq \frac{C}{N} \sum_{1 \leq k \leq k_1 } \frac{N^{(1)}_k}{ 2^k G} + C \notag \\
&\leq C k_1 + C \leq C \log (N).
\end{align}
For the term $A_2$ we define the intervals for $k \in \zz$
\beq
\Itwo_k := E + T \Re [ \mfcl (z) ] + 2 Tk + [-T, T]
\eeq
and let 
\beq
N^{(2)}_k := | \{ i : V_i \in \Itwo_k \} |.
\eeq

Let $k_2$ and $k_3$ be the index of the leftmost and rightmost intervals $\Itwo_k$ that intersect $\I_{E_0, \sqrt{q}G}$, respectively.
By Lemma \ref{lem:vbound}, 
\beq
N^{(2)}_k \leq CTN, \qquad \mbox{for }k_2 \leq k \leq k_3.
\eeq
Then, 
\begin{align}
A_1 &\leq \frac{1}{N} \sum_{k_2 \leq k \leq k_3} \sum_{V_i \in \Itwo_k } |g_i (z) | \notag \\
&\leq \frac{C}{N} \sum_{k_2 \leq k \leq k_3} \frac{ N^{(2)}_k }{ T ( |k| +1 )} \leq C \log N.
\end{align}
\qed

We also require some results about the qualitative properties of the time evolution of $\mfct$.

\bel \label{lem:timeder} 
Let $\mfct$ satisfy (\ref{eqn:mfct}).  Then
\beq
\del_t \mfct (z) = \frac{1}{2} \del_z [ \mfct (z) ( \mfct (z) + z ) ]. \label{eqn:complexburger}
\eeq
We also have
\beq
\del_t \gamit = - \Re [ \mfct ( \gamit ) ] - \frac{1}{2} \gamit \label{eqn:gamtime}.
\eeq
Fix $\om > 0 $ and $q >0$.  We have uniformly for $i \in \Amuto$ and $t, t_0 \in \Tom$,
\beq
|\del_t \gamit | \leq C \log N. \label{eqn:gambd}
\eeq
as well as
\beq
| \del_t \rhofct ( E ) | \leq \frac{C}{t} \label{eqn:rhotbd}
\eeq
for $E \in \I_{E_0, q G}$.

\eel
\proof Equation (\ref{eqn:complexburger}) is a straight-forward exercise differentiating (\ref{eqn:mfct}).  The equation (\ref{eqn:gamtime}) follows as in the proof of Lemma 4.4 of \cite{Kevin2}. The resulting bound (\ref{eqn:gambd}) follows from the bound $| \mfct | \leq C \log N$ which is obtained from Lemma \ref{lem:gilogbd}.  The bound (\ref{eqn:rhotbd}) follows from (\ref{eqn:complexburger}) and the fact that we bounded $| \del_z \mfct|  \leq C / t$ in the proof of Lemma \ref{lem:rstab}. \qed

\noindent{\bf Proof of Lemma \ref{lem:bulkindices}. } This follows directly from (\ref{eqn:gambd}).  
\qed

\subsection{Local deformed semicircle law} 
\subsubsection{Preliminaries}

We collect here some identities and tools required for our proofs.  
For a set of indices $\tt \subset [1, N]$ we denote the minor of $H_T$ obtained by removing the columns and rows in $\tt$ by $H_T^\btt$.  Similarly we denote the Green's function of $H_T^\btt$ by $G^\btt_{ij}$. For the minors we will still normalize the empirical Stieltjes transform as
\beq
\mN^\btt (z) = \frac{1}{N} \tr \frac{1}{ H_T^\btt -z }
\eeq  
We use the notation $\tt i  = \tt \cup \{ i \}$ and $H_T^{(i)}$, $G^{(i)}_{jk}$ for $\tt = \{ i \}$.  We use the notation
\beq
\sum^{\btt}_{ij} := \sum_{\substack{i \notin \tt \\ j \notin \tt}}
\eeq
to denote summation over all indices not in $\tt$.

For a row $i$ and random variable $X$ we denote $P_i [X] = \ee [ X | H_T^{(i)} ]$ and $Q_i [X] = X - P_i [X]$.  

In the following lemma we collect the resolvent identities required for the proof.  We refer the reader to, for example, \cite{renyione} for a proof.
\bel 
The Schur complement formula is, for $i \notin \tt$:
\beq \label{eqn:schur1}
G_{ii}^{\btt} = \frac{1}{ h_{ii}- z - \sum_{k, l}^{\btti} h_{ik} G_{kl}^{\btti} h_{li} }.
\eeq
For $i, j, k \notin \tt$ and $i, j \neq k$ we have
\beq \label{eqn:res1}
G^\btt_{ij} = G_{ij}^{( \tt k )} + \frac{ G_{ik}^\btt G_{kj}^\btt}{G_{kk}^\btt} , \quad \frac{1}{G_{ii}^\btt} = \frac{1}{G_{ii}^{\bttk}} - \frac{ G_{ik}^\btt G_{ki}^\btt}{ G_{ii}^\btt G_{ii}^\bttk G_{kk}^\btt }.
\eeq
For $i, j \notin \tt$ and $i \neq j$, 
\beq \label{eqn:res2}
G_{ij}^\btt = - G_{ii}^\btt \sum_k^\btti h_{ik} G_{kj}^\btti = - G_{jj}^\btt \sum_k^\bttj G_{ik}^{\bttj} h_{kj}
\eeq
and
\beq \label{eqn:res3}
G_{ij}^\btt = - G_{ii}^\btt G_{jj}^\btti ( h_{ij} - \sum_{m, n}^{ ( \tt i j ) } h_{im} G_{mn}^{ ( \tt i j ) } h_{nj } ).
\eeq
We have also the Ward identity
\beq \label{eqn:ward}
\sum_{j} | G_{ij}^\btt |^2 = \frac{1}{ \eta} \Im G_{ii}^\btt.
\eeq

\eel

We collect here the large deviations estimates we will require.  We refer the reader to, for example, \cite{erdos2012bulk} for a proof.
\bel Let $(a_i)$ and $(b_i)$ be centered independent random variables with variance $\sigma^2$ and having subexponential decay
\beq
\pp ( |a_i | \geq x \sigma ) \leq C_0 \e^{- x^{ 1 / \theta } } , \qquad\pp ( |b_i | \geq x \sigma ) \leq C_0 \e^{- x^{ 1 / \theta } }
\eeq
for some positive constant $C_0$ and $ \theta >1$.  Let $A_i \in \cc$ and $B_{ij} \in \cc$.  Then there exist constants $a_0 >1$, $A_0 \geq 10$ and $C\geq 1$ depending on $\theta$ and $C_0$ s.t. for $a_0 \leq \xi \leq A_0 \log \log N$ and $\varphi = ( \log N)^C$, 
\begin{align}
\pp \left( \left| \sum_{i=1}^N A_i a_i \right| \geq \varphi^\xi \sigma \left( \sum_{i=1}^N |A_i|^2 \right)^{1/2}\right) &\leq \e^{ - ( \log N)^\xi } \label{eqn:ld1} \\
\pp \left( \left| \sum_{i=1}^N \bar{a}_i B_{ii} a_i - \sum_{i=1}^N \sigma^2 B_{ii} \right| \geq \varphi^\xi \sigma \left( \sum_{i=1}^N |B_{ii}|^2 \right)^{1/2}\right) &\leq \e^{ - ( \log N)^\xi } \label{eqn:ld2} \\
\pp \left( \left| \sum_{i \neq j}^N \bar{a}_i B_{ij} a_j \right| \geq \varphi^\xi \sigma \left( \sum_{i \neq j}^N |B_{ij}|^2 \right)^{1/2}\right) &\leq \e^{ - ( \log N)^\xi }  \label{eqn:ld3} \\
\pp \left( \left| \sum_{i, j=1}^N \bar{a}_i B_{ij} b_j \right| \geq \varphi^\xi \sigma \left( \sum_{i=1}^N | B_{ij}|^2 \right)^{1/2}\right) &\leq \e^{ - ( \log N)^\xi } \label{eqn:ld4}
\end{align}
for $N$ sufficiently large.
\eel

\subsubsection{Proof of Theorem \ref{thm:ll}}

 The proof follows closely the proof of Theorem 2.10 in \cite{Kevin1}.  We will just note the differences here.  The first difference is that we consider $V$ fixed instead of random.  This change was already dealt with in the proof of Theorem 3.3 of \cite{Kevin2}.  The extra randomness in $V$ corresponds to the larger error terms which appear throughout \cite{Kevin1}.  For us, taking $V$ fixed corresponds to setting $\lambda =0$ everywhere in the estimates appearing in Sections 3 and 4 of \cite{Kevin1}.  Additionally, in the paper \cite{Kevin1} the local law is also proved at the edge of the spectrum.  We will not deal with the edge case and will therefore avoid some of the technicalities encountered in \cite{Kevin1}.  We split the proof of the local law into two parts, first for $z \in \D_1$ and then for $z \in \D_2$.  The proof in the domain $\D_1$ is similar to \cite{Kevin1} and we now summarize the changes.
 
 \subsubsection{Local law in the spectral domain $\D_1$}

The major difference between the proof here and that in \cite{Kevin1} is of course that we consider $H_T = V + \sqT W$ for small $T$ whereas in \cite{Kevin1}, $T$ is order $1$.  We summarize what difficulties this causes and how they are dealt with.  The first lies in the stability estimate (\ref{eqn:stab2})  which only gives a lower bound of $c T$ whereas the analogue in \cite{Kevin1} is Corollary A.2 where the lower bound is $c$.  This changes a few steps of the proof.  Instead of (3.17) of \cite{Kevin1}, the Schur complement formula is written as
\beq \label{eqn:schur}
G_{ii} = \frac{1}{ V_i - z - T \mfcl + T \psi + Y_i }
\eeq
where
\beq
\psi = \mfcl - \mN , \qquad Y_i = \sqT w_{ii} + T ( \mN - \mN^{(i)} ) - Z_i
\eeq
and
\beq
Z_i = \sum_{k, l}^{(i)} h_{ik} G_{kl}^{(i)} h_{li} - T \mN^{(i)} = Q_i \left[ \sum_{k, l}^{(i)} h_{ik} G_{kl}^{(i)} h_{li}  \right]
\eeq
with $\mN^{(i)}$ the empirical Stieltjes transform of the $i$th minor of $H_T$.

Let
\beq
\Lam (z) := | \psi | = | \mN (z) - \mfcl (z) |, \qquad \Psi (z) := \varphi^{\xi} \sqrt{ \frac{ \Lam (z) + 1}{ N \eta } }.
\eeq
Define
\beq
g_i (z) := \frac{1}{ V_i - z -T \mfc (z) }.
\eeq

We have the following a priori bounds on $Z_i$ and $Y_i$ on the event that $\Lam (z) = o (1)$.  This is our version of Lemma 3.8 of \cite{Kevin1}.  

\bel \label{lem:apriori} On the event $\{ \Lam (z) \leq \varphi^{ - 2 \xi } \}$ we have with $(\xi, \nu)$-high probability,
\beq
\max_i |Z_i| \leq C T \Psi, \qquad \max_i | g_i Y_i | \leq C \Psi \label{eqn:apriori}
\eeq
and
\beq
\max_i |Y_i| \leq C \left( \frac{\varphi^{\xi} \sqT}{ \sqrt{N} } + T \Psi \right) \label{eqn:yapriori}
\eeq
for a constant $C$.
\eel
\proof By the large deviations bounds (\ref{eqn:ld2}), (\ref{eqn:ld3}) we have with $(\xi, \nu)$-high probability
\begin{align}
|Z_i| \leq C \varphi^\xi T \left( \frac{1}{N^2} \sum_{k, l}^{(i)} | G^{(i)}_{kl} |^2 \right)^{1/2} \leq C T \Psi
\end{align}
where we have used (\ref{eqn:ward}) and the inequality \cite{Kevin1}
\beq
|\mN - \mN^{(i)} | \leq (N \eta)^{-1}. \label{eqn:mni}
\eeq
  The bound (\ref{eqn:yapriori}) for $Y_i$ follows from this and the subexponential decay assumption for $w_{ii}$ which yields with $(\xi, \nu)$-high probability,
\beq
\sqT |w_{ii}| \leq \sqT \frac{\varphi^{\xi}}{ \sqrt{N}} . \label{eqn:wsubexp}
\eeq
For the bound (\ref{eqn:apriori}) for $Y_i$ we use (\ref{eqn:mni}), (\ref{eqn:wsubexp}) and $|g_i | \leq \min \{ \eta^{-1} , T^{-1} \}$. \qed

Using these estimates, we can prove the following analog of the weak-self consistent equation which is Lemma 3.9 of \cite{Kevin1}.
\bel \label{lem:ws} On the event $\{ \Lam (z) \leq \varphi^{- 2 \xi } \}$ we have with $(\xi, \nu)$-high probability,
\beq
\left| ( 1 - T R_2 (z) ) \psi - T^2 R_3 (z) \psi^2 \right| \leq C \Psi + C \frac{ \Lam ^2}{ \log N }. \label{eqn:weakself}
\eeq
\eel
\proof By Lemma \ref{lem:apriori} we have on the event $\{ \Lam (z) \leq \varphi^{-2 \xi } \}$ that $|Y_i| \leq T o (1)$. Therefore, by the stability bounds (\ref{eqn:stab2}) and (\ref{eqn:dumbstab}) we can expand (\ref{eqn:schur}) and obtain
\begin{align}
\frac{1}{N} \sum_i G_{ii} &= \frac{1}{N} \sum_i g_i + \frac{1}{N} \sum_i g_i^2 ( T \psi + Y_i ) + \frac{1}{N} \sum_i g_i^3 T^2 \psi^2 \notag \\
&+ \frac{1}{N} \sum_i g_i^3 ( 2 T \psi Y_i + Y_i^2 ) + \O \left( ( \max_i |g_i|^3 ( |Y_i|^3 + T^3 \Lam^3 ) ) N^{-1} \sum_j |g_j | \right). \label{eqn:ws1}
\end{align}
 We use a simple Schwarz inequality $| \psi Y_i| \leq C |Y_i |^2 (\log N) + \log (N)^{-1} \Lam^2$  on the first term in the second line of (\ref{eqn:ws1}).  By Lemmas \ref{lem:apriori} and \ref{lem:gilogbd} we have
\beq
\frac{1}{N} \sum_i |g_i |^{1+r} |Y_i |^r \leq  C \log (N) \Psi^{r}, 
\eeq
for $r \geq 0$.  This, together with (\ref{eqn:apriori}) yields the claim. \qed

With this, one can derive the weak local deformed semicircle law. 
\bet[Weak local deformed semicircle law] The event
\beq
 \bigcap_{ T \in \Tom, z \in \D_1 } | \mN (z) - \mfcl (z) | \leq C \frac{ \varphi^{\xi}}{ ( N \eta )^{1/3}}
\eeq
holds with $(\xi, \nu)$-high probability.

\eet
This is derived from the weak self-consistent equation exactly as Theorem 3.1 is derived in \cite{Kevin1}. To be more precise, one proves the initial estimates for $\eta = 2$ in the same fashion as Lemma 3.10 of \cite{Kevin1}.  The dichotomy argument of Lemma 3.12 is easier as we have that the LHS of (3.45) of \cite{Kevin1} can be taken to be $1$, and so one only needs to do the `Bulk case' argument in the proof of Lemma 3.12 of \cite{Kevin1}; i.e., we will only be considering the regime $\eta \geq \tilde{\eta} ( U, E)$. The stability of the coefficients in the self-consistent equation is provided by the estimates of Lemmas \ref{lem:rstab} and \ref{lem:dumbstab} (as we are in the bulk the coefficient $ ( 1  - T R_2 (z) )$ does not degenerate which means that $\alpha \geq c$ where $\alpha$ is defined as in \cite{Kevin1}).  

The proof of the weak local deformed semicircle law is then completed as in \cite{Kevin1}.  We record here the following corollary of the above proof.  We will not need it in this paper but record it for possible future use.

\bec
 Suppose that one has that $c \leq \Im [ m_V (E + \i \eta ) ] \leq C$ for $\eta \geq \eta_0$ where $\eta_0$ satisfies (for example) $N^{-1} \leq \eta_0 \leq N^{-c}$, and all $E \in \I$ where $\I$ is some fixed interval.  Let $\eta_j = \varphi^{10 j \xi} \eta_0$ for $j=1, 2$.  Then with $(\xi, \nu)$-high probability, one has
\begin{align}
| m_N (z) - \mfcll (z) | \leq C \frac{ \varphi^\xi}{ ( N \eta )^{1/3} }
\end{align}
for all $0 \leq T \leq N^{-c}$, $E \in \I'$ and $\eta_2 \leq \eta \leq 10$ where $\I'$ is an interval compactly contained in $\I$.

\eec
\proof  It is easy to see that the above proof directly implies the weak law for $T \geq \eta_1$ with $T \ll G^2$ for $\eta$ down to the optimal scale $\eta \gtrsim N^{-1}$.   On the other hand, the above proof also yields the weak law for $0\leq T \leq \eta_1$ for $\eta \geq \eta_2$.  For this, the key input is that the bounds (\ref{eqn:stab1}), (\ref{eqn:stab2}), and (\ref{eqn:stab4}) hold in the regime $\eta \geq \eta_2$.  The first two are obvious as we have $\Im [ \mfcl ] = \Im [ m_V ( z +T \mfcl ) ]$ and $\Im [ z + T \mfcl ] \geq \eta_2 \gg \eta_0$.  The upper bounds of (\ref{eqn:stab4}) are clear as $\eta_2 \gg T$.  The lower bound is slightly nontrivial but we have,
\beq
\Re [ 1 - T R_2 (z) ] \geq 1 - \frac{T}{ \eta_2} \Im [ m_V ( z + T \mfcl (z) ) ] \geq \frac{1}{2}
\eeq
for large enough $N$.   The expansion in the proof of Lemma \ref{lem:ws} is still possible as now one has the stability bound $|v_i -z - T \mfc (z) | \geq \eta_2$ and $|T \psi| + |Y_i| \leq T o(1) + (T/N)^{1/2} \varphi^{\xi} \leq \eta_2 o(1)$.

\qed

We now return to the proof of Theorem \ref{thm:ll}.  In order to get the strong self-consistent equation, Lemma \ref{lem:sself} below, we require the following fluctuation averaging lemma.  As the proof is somewhat lengthy, we defer it to the next subsection.  We introduce deterministic control parameters $\gamma (z)$ and $\Phi (z)$ by requiring
\beq
\gamma (z) \leq \varphi^{- 2 \xi }, \qquad \Phi (z) = \varphi^{\xi} \sqrt{ \frac{ 1 + \gamma (z) }{ N \eta}}.
\eeq

\bel \label{lem:fa} Suppose that the event
\beq
\Lam (z) \leq \gamma (z), \quad \forall z \in \D_{1}
\eeq
holds with $(\xi, \nu)$-high probability.  Then,
\beq
\left| \frac{1}{N} \sum_{i=1}^N g_i^2  Q_i \left( \frac{1}{G_{ii}} \right)   \right| \leq C \varphi^{ 10 \xi } \Phi^2.
\eeq
with $(\xi -2, \nu)$-high probability.

\eel

The above lemma allows us to deduce the strong self-consistent equation.  It is our version of Lemma 4.5 of \cite{Kevin1}.
\bel \label{lem:sself} Suppose that the event
\beq
\Lam (z) \leq \gamma (z), \quad \forall z \in \D_{1}
\eeq
holds with $(\xi, \nu)$-high probability.  Then,
\beq
| ( 1 - T R_2 ) \psi -T^2 R_3 \psi^2 | \leq C \frac{ \Lam^2}{ \log N} + C \varphi^{ 10 \xi } \Phi^2
\eeq
with $(\xi-2, \nu)$-high probability.

\eel

\proof We proceed as in the proof of Lemma \ref{lem:ws} and expand (\ref{eqn:schur}).  We get
\begin{align}
\frac{1}{N} \sum_i G_{ii} &= \frac{1}{N} \sum_i g_i + \frac{1}{N} \sum_i g_i^2 ( T \psi + Y_i ) + \frac{1}{N} \sum_i g_i^3 T^2 \psi^2 \notag \\
&+ \frac{1}{N} \sum_i g_i^3 ( 2 T \psi Y_i + Y_i^2 ) + \O \left( ( \max_i |g_i|^3 ( |Y_i|^3 + T^3 \Lam^3 ) ) N^{-1} \sum_j |g_j | \right). \label{eqn:ss1}
\end{align}
Note that we have
\beq
Q_i \left[ \frac{1}{G_{ii} } \right] = Q_i \left[ v_i + \sqT w_{ii} - z - \sum_{k, l}^{(i)} h_{ik} G_{kl}^{(i)} h_{li} \right] = \sqT w_{ii} - Z_i. \label{eqn:ab1}
\eeq
Hence,
\beq
Y_i = T ( \mN - \mN^{(i)} ) + Q_i \left[ \frac{1}{G_{ii} } \right] .
\eeq
Therefore, by Lemma \ref{lem:fa} and (\ref{eqn:mni}) we have
\beq
\left| \frac{1}{N} \sum_i g_i^2 Y_i \right| \leq  C\varphi^{10 \xi }  \Phi^2
\eeq
with $(\xi-2, \nu)$-high probability.  The other terms on (\ref{eqn:ss1}) are bounded as in the proof of Lemma \ref{lem:ws} using the a-priori bounds in Lemma \ref{lem:apriori}. \qed

The proof of Theorem \ref{thm:ll} now follows from the strong self-consistent equation as it does in the proof of Theorem 2.10 at the end of Section 4 in \cite{Kevin1}.

\subsubsection{Proof of Lemma \ref{lem:fa}}

Lemma \ref{lem:fa} follows directly from the following moment bound and Chebyshev's inequality.
\bel If the event \label{lem:momentbd}
\beq
\Lam (z) \leq \gamma (z), \quad \forall z \in \D_{1}
\eeq
holds with $(\xi, \nu)$-high probability, then
\beq
\ee \left| \frac{1}{N} \sum_{i} g_i^2 Q_i \left( \frac{1}{G_{ii} } \right) \right|^{2r} \leq (C r )^{C r } ( \varphi^{ 3 \xi } \Phi )^{4r}
\eeq
for every $r \leq \nu ( \log N )^{ \xi - 3/2 } /2$.

\eel

For its proof we will require the following a-priori bounds on the Green's function matrix elements.
\bel \label{lem:apriori2}
Suppose that the event
\beq
\Lam (z) \leq \gamma (z), \quad \forall z \in \D_{1} 
\eeq
holds with $(\xi, \nu)$-high probability.  We then have uniformly for $| \tt | \leq ( \log N )^{ \xi -1 }$,
\beq
\left| \frac{g_i}{G_{ii}^{\btt} } -1  \right| \leq \varphi^{  \xi } \Phi  + \varphi^{- \xi}\label{eqn:ap1}
\eeq
and for $0 \leq s \leq 1$ and $i \neq j$
\beq
\frac{ |G_{ij}^{\btt} |}{ |g_i|^s |g_j|^{1-s} } \leq \varphi^{ 2\xi  } \Phi \label{eqn:ap2}
\eeq
and
\beq
|g_i| \left| Q_i \left[ \frac{1}{ G_{ii}^{\btt} } \right] \right| \leq \varphi^{  2\xi } \Phi. \label{eqn:ap3}
\eeq

\eel
\proof From the Schur complement formula we have
\beq
\left| \frac{g_i}{ G_{ii}^{\btt} } -1 \right| = |g_i| \left| T \psi + Y^{\btt}_i \right|
\eeq
where
\beq
Y^{\btt}_i := \sqT w_{ii} + T ( \mN - \mN^{ \btti} ) - Z_i^{\btt} , \quad Z_i^{\btt} := \sum_{k, l}^{\btti} h_{ik} G_{kl}^{\btti} h_{li} - T \mN^{\btti} .
\eeq
The proof of Lemma \ref{lem:apriori} extends to $Y_i^{\btt}$ and $Z_i^{\btt}$ uniformly for $| \tt | \leq ( \log N )^{ \xi-1}$ and we conclude (\ref{eqn:ap1}).  Moreover, we see that this observation also proves (\ref{eqn:ap3}) after using
\beq
Q_i \left[ \frac{1}{G_{ii}^{\btt} } \right] = \sqT w_{ii} - Z_i^{\btt}.
\eeq

For (\ref{eqn:ap2}) we start with (\ref{eqn:res3}) and obtain with $(\xi, \nu)$-high probability by the large deviations bounds, (\ref{eqn:ap1}) and (\ref{eqn:ward})
\begin{align}
| G_{ij}^{\btt} | &= | G_{ii}^{\btt} G_{jj}^{ \btti } | \left|  h_{ij} - \sum_{m, n}^{  \bttij  } h_{im} G_{mn}^{\bttij} h_{nj}  \right| \notag \\
&\leq \varphi^{\xi} |g_i g_j | \left( \frac{\sqT}{ \sqrt{N}} + T \sqrt{ \frac{ \Im \mN^{ \bttij}  }{ N \eta } } \right) \notag \\
&\leq \varphi^{2 \xi } |g_i|^s |g_j |^{1-s} \Phi.
\end{align}
In the last line we have used $|g_l| \leq T^{-1}$ or $|g_l| \leq \eta^{-1}$ as appropriate.  \qed

\noindent{\bf Proof of Lemma \ref{lem:momentbd}.}   We illustrate the method by doing $r=1$. We have,
\begin{align}
\ee \left| \frac{1}{N} \sum_i g_i^2 Q_i \left[ G_{ii}^{-1} \right] \right|^2 &= \frac{1}{N^2} \sum_{i, j} \ee \left[\bar{g}_i^2 Q_i \left[ \bar{G}_{ii}^{-1} \right] g_j^2 Q_j \left[ G_{jj}^{-1} \right] \right] \notag \\
&= \frac{1}{N^2} \sum_{i=1}^N |g_i|^4 \ee \left| Q_i \left[ G_{ii}^{-1} \right] \right|^2 + \frac{1}{N^2} \sum_{i \neq j} \ee \left[\bar{g}_i^2 Q_i \left[ \bar{G}_{ii}^{-1} \right] g_j^2 Q_j \left[ G_{jj}^{-1} \right] \right] \notag \\
&= A_1 + A_2.
\end{align}
For the first sum we apply (\ref{eqn:ap3}). Hence,
\beq
A_1  \leq \varphi^{ 4 \xi } \Phi^2 \frac{1}{N^2} \sum_{i=1}^N |g_i|^2 \leq \varphi^{ 4 \xi } \Phi^2 \frac{1}{ N \eta } \leq ( \varphi^{ 2 \xi } \Phi )^4
\eeq
where we have used
\beq
\frac{1}{N} \sum_{i=1}^N |g_i|^2 \leq  \frac{1}{ \eta}.
\eeq
We apply the second identity of  (\ref{eqn:res1}) to the summand in $A_2$.  We have
\begin{align}
\ee \left[\bar{g}_i^2 Q_i \left[ \bar{G}_{ii}^{-1} \right] g_j^2 Q_j \left[ G_{jj}^{-1} \right] \right] &=\ee \left[\bar{g}_i^2 Q_i \left[  \frac{1}{\bar{G}_{ii}^{(j)} } - \frac{ \bG_{ij} \bG_{ji} }{ \bG_{ii} \bG_{ii}^{(j)} \bG_{jj} }  \right] g_j^2 Q_j \left[ \frac{1}{G_{jj}^{(i)} } - \frac{ G_{ji} G_{ij} }{ G_{jj} G_{jj}^{(i)} G_{ii} }\right]  \right] \notag \\
&= \ee \left[\bar{g}_i^2 Q_i \left[ \frac{ \bG_{ij} \bG_{ji} }{ \bG_{ii} \bG_{ii}^{(j)} \bG_{jj} }  \right] g_j^2 Q_j \left[  \frac{ G_{ji} G_{ij} }{ G_{jj} G_{jj}^{(i)} G_{ii} }\right]  \right] 
\end{align}
where we have used
\beq
\ee [ X Q_i [ Y ] ] = \ee [ Q_i [ X Y ] ] = 0 
\eeq
for $X$ measureable wrt to the $i$th minor of $H_T$.  From Lemma \ref{lem:apriori2} we see that 
\beq \label{eqn:resbds}
\left| \frac{g_l}{ G_{ll}^{ ( \tt_1 ) } } \right| \leq 2 , \qquad  \left| \frac{ G_{lk}^{ ( \tt_1 ) }}{ G_{kk}^{  ( \tt_2 ) } } \right| \leq \varphi^{ 2 \xi } \Phi, \qquad \left| \frac{ G_{lk}^{ ( \tt_1 ) }}{ G_{ll}^{  ( \tt_2 ) } } \right| \leq \varphi^{ 2 \xi } \Phi , \qquad l \neq k
\eeq
holds with $(\xi-2, \nu)$-high probability uniformly for $| \tt_1|, | \tt_2 | \leq ( \log N )^{ \xi -1 }$.  From these inequalities we deduce 
\beq
\left| \ee \left[\bar{g}_i^2 Q_i \left[ \bar{G}_{ii}^{-1} \right] g_j^2 Q_j \left[ G_{jj}^{-1} \right] \right]  \right| \leq |g_i g_j | ( \varphi^{ 2 \xi } \Phi )^4
\eeq
and so
\beq
A_2 \leq ( \varphi^{2 \xi } \Phi )^4 \frac{1}{N^2} \sum_{i, j} |g_i g_j | \leq ( \varphi^{3 \xi } \Phi )^4
\eeq
where we have used the fact from Lemma \ref{lem:gilogbd} that
\beq
\frac{1}{N} \sum_i |g_i | \leq C \log N.
\eeq
This proves the case $r=1$.  

For the general case we will just explain how to adapt the proof of (4.5) of \cite{Kevin1} to our set-up.  We start with (4.28) of \cite{Kevin1}:
\beq \label{eqn:bigsum}
\ee \left| \frac{1}{N} \sum_{i=1}^N g_i^2 Q_i \left( \frac{1}{ G_{ii} } \right) \right|^{2r} = \frac{1}{N^{2r} } \sum_{ \Gamma \in \P_{2r} } \sum_{i_1, ..., i_r} \1_{ \{ \Gamma = \Gamma ( \underline{i} ) \} } \ee [ g^2_{i_1} Q_{i_1} ( \widebar{ G_{i_1, i_1} ^{-1} } ) ... g^2_{i_r} Q_{i_{2r}} (  G_{i_r, i_r} ^{-1}  ) ]
\eeq
where $\P_{2r}$ denotes the set of partitions on $2r$ letters and $\Gamma (\underline{i} )$ denotes the element of $\P_{2r}$ defined by the equivalence relation $a \sim b$ iff $i_a = i_b$, and $\underline{i} = (i_1, ..., i_{2r} )$.

Fix now a partition $\Gamma$ and let $\underline{i}$ satisfy $  \Gamma ( \underline{i}) = \Gamma$. We apply the same algorithm as in \cite{Kevin1} to the summand on the RHS of (\ref{eqn:bigsum}).  We omit the precise details of the algorithm.  The result is the expansion into the sum
\beq
\ee [ g^2_{i_1} Q_{i_1} ( \widebar{ G_{i_1, i_1} ^{-1} } ) ... g^2_{i_r} Q_{i_{2r}} (  G_{i_r, i_r} ^{-1}  ) ] = \sum_{\sigma_1, ..., \sigma_{2r} } \ee \left[g_{i_1}^2 Q_{i_1} \widebar{ ( F_{i_1} )_{\sigma_1}  } ... g_{i_{2r} }^2 Q_{i_{2r}} ( F_{i_{2r} } )_{ \sigma_{2r} }   \right] \label{eqn:ab2}
\eeq
where the $\sigma_k$ run over the finite binary sequences generated by the algorithm and $ (F_{i_k} )_{\sigma_k }$ are the associated monomials in the resolvent entries.  We will eventually apply (\ref{eqn:resbds}) to the monomials $\Fk$ in the above expression.  To see that we can do so, we require the following lemma.  It is an easy consequence of the definition of the algorithm generating the monomials which uses the identities (\ref{eqn:res1}). We give the full details of the proof.  We encourage the reader unfamiliar with the algorithm to consult \cite{Kevin1} or Appendix B of \cite{general}.

\bel \label{lem:gijmatch} Let $( F_{i_k} )_{\sigma_k}$ be a monomial generated by the algorithm outlined after (4.29) in \cite{Kevin1}.  Let $\F_1 $ be the set of off-diagonal resolvent entries in the numerator and $\F_2 $ be the set of diagonal resolvent entries in the denominator.  Suppose that $\F_1 \neq \emptyset$.  Note that $|\F_2| = | \F_1 | + 1$. There exists an injective function $\pi : \F_1 \to \F_2$ such that,
\begin{enumerate}
\item For $G_{xy}^{\btt} \in \F_1$, the lower index of $\pi ( G_{xy}^{\btt})$ is either $x$ or $y$.
\item The unique diagonal resolvent entry $G^{\btt}_{xx} \in \F_2 \backslash \pi ( \F_1)$ not lying in the image of $\pi$ satisfies $x = i_k$.  That is, its lower index is the same as the index $Q_{i_k}$ over which the partial expectation is taken.

\end{enumerate}

\eel
\proof Recall that the off-diagonal resolvent entries are generated by applying either of the rules in (\ref{eqn:res1}) to an existing monomial. Let $\sigma_k = a_1 a_2 ... a_n$ for $a_i \in \{ 0, 1 \}$.  For each step $j$ we construct $\pi_j : \F_{1, j} \to \F_{2, j}$ satisfying the lemma, where $\F_{1, j}$ and $\F_{2, j}$ are the off-diagonal and diagonal resolvent entries of $( F_{i_k} )_{a_1 ... a_j }$, respectively.

    If the first step $a_1 = 0$ then there is nothing to do as there are no off-diagonal resolvent entries yet. If $a_1 = 1$ then 
\beq
( F_{i_k } )_{a_1} = - \frac{ G_{i_k j} G_{j, i_k } } { G_{i_k i_k } G_{i_k i_k }^ {(j)} G_{jj} }
\eeq
for some index $j$. We can take
\beq
\pi_1 ( G_{i_k j}  ) = G_{i_k i_k }^ {(j)} , \qquad \pi_1 ( G_{j, i_k } ) = G_{jj}.
\eeq
We proceed by induction. Assume that we have constructed $\pi_j$.  First consider the case $a_{j+1} = 0$. If the second rule of (\ref{eqn:res1}) was applied then we can let $\pi_{j+1} = \pi_j$ if the rule was applied to a diagonal resolvent entry not lying in the image of $\pi_j$.  If the rule substituted $G_{ii}^{(\tt_1)} \to G_{ii}^{(\tt_1 l)}$ with $\pi_j ( G_{xy}^{(\tt_2)} )= G_{ii}^{(\tt_1)}$, then we only have to define $\pi_{j+1} ( G_{xy}^{(\tt_2)} ) = G_{ii}^{( \tt_1 l)}$.  For all other off-diagonal entries we just let $\pi_{j+1} = \pi_j$, because the other elements of $\F_{1, j+1}$ coincide with the other elements of $\F_{1, j}$ and the other diagonal resolvent entries remain unchanged.

  On the other hand suppose that the first rule was applied to $G_{xy}^{(\tt)}$ which became $G_{xy}^{(\tt l ) }$.  We can let $\pi_{j+1} ( G_{xy}^{ ( \tt l )} ) = \pi_j ( G_{xy}^{( \tt )})$.  The other elements of $\F_{1, j+1}$ coincide with the other elements of $\F_{1, j}$ and we set $\pi_{j+1} = \pi_j$ on these elements.

Now we can suppose that $a_{j+1} =1$.  Suppose that the first rule was applied and that we have the replacement
\beq
G_{xy}^{ \btt } \to \frac{ G_{xl}^{ ( \tt l ) } G_{ l y}^{ ( \tt l ) } } { G_{ll}^{ \btt } }.
\eeq
Suppose that the lower index of $\pi_j ( G_{xy}^{ \btt } )$ was $x$.  Then we set,
\beq
\pi_{j+1} ( G_{xl}^{ ( \tt l ) } ) = \pi_j ( G_{xy}^{ \btt } ) , \qquad \pi_{j+1} ( G_{ l y}^{ ( \tt l ) }  ) = G_{ll}^{ \btt } 
\eeq
and set $\pi_{j+1} = \pi_j$ on the other elements of $\F_{1, j+1}$ which coincide with the other elements of $\F_{1, j}$.  The case when the lower index of $\pi_j ( G_{xy}^{ \btt } )$ is $y$ is similar.

If the second rule was applied then we have the replacement
\beq
\frac{1}{ G_{xx}^{\btt}} \to - \frac{ G_{x y}^{ \btt}  G_{y x }^{ \btt } } { G_{xx}^{ \btt} G_{xx}^{ ( \tt y ) } G_{yy}^{ \btt } }.
\eeq
We then set
\beq
\pi_{j+1} ( G_{x y}^{ \btt} ) = G_{xx}^{ ( \tt y ) } , \qquad  \pi_{j+1} (  G_{y x }^{ \btt } ) = G_{yy}^{ \btt }
\eeq
and set $\pi_{j+1} = \pi_j$ on the other elements of $\F_{1, j+1}$. This yields the lemma. \qed

As a corollary we have
\beq \label{eqn:onesbd}
\left| g_{i_k} Q_{i_k} ( F_{i_k} )_{\sigma_k } \right| \leq  ( \varphi^{2 \xi } \Phi )^{ b (\sigma_k ) + 1}
\eeq
where $b (\sigma_k )$ is the number of ones in the string $\sigma_k$. If $b ( \sigma_k )  = 0$ this follows from (\ref{eqn:ap3}). If $b (\sigma_k ) > 0$ then note that the number of off-diagonal entries in the numerator equals $b ( \sigma_k ) +1$.  The bound (\ref{eqn:onesbd}) then follows from Lemma \ref{lem:gijmatch} and (\ref{eqn:resbds}).
 
   For a label $a \in \{1, ..., 2r \}$ we let $[a]$ denote the block of $a$ in $\Gamma$.  We let $S ( \Gamma ) := \{ a : | [ a ] | = 1 \}$ denote the set of single labels in $\Gamma$ and denote by $s = | S ( \Gamma ) |$ its cardinality.  We wish to prove that any nonzero term on the RHS of  (\ref{eqn:ab2}) satisfies 
\beq \label{eqn:sbd}
\left| g_{i_1} Q_{i_1} \widebar{ ( F_{i_1} )_{\sigma_1}  } ... g_{i_{2r}}Q_{i_{2r}} ( F_{i_{2r} } )_{ \sigma_{2r} }  \right| \leq C ( \varphi^{2 \xi } \Phi )^{2r + s }.
\eeq
First, suppose that one of the monomials, say $(F_{i_k})_{\sigma_k}$, is not maximally expanded (recall that a monomial is maximally expanded if each resolvent $G_{xy}^{\btt}$ in the monomial satisfies $S ( \Gamma ) \subseteq \{ x, y, \tt \}$).  By the definition of the algorithm generating the monomials, we must have that $\Fk$ contains at least $4r$ off-diagonal resolvent entries in the numerator.  Then $b ( \sigma_k ) \geq 4r -1$ and so from (\ref{eqn:onesbd}) we have
\beq
\left| g_{i_1} Q_{i_1} \widebar{ ( F_{i_1} )_{\sigma_1}  } ... g_{i_{2r}} Q_{i_{2r}} ( F_{i_{2r} } )_{ \sigma_{2r} } \right| \leq ( \varphi^{2 \xi } \Phi )^{1 + 1 + .... + (4r ) + 1 ...} \leq ( \varphi^{2 \xi } \Phi )^{4 r} \leq  ( \varphi^{2 \xi } \Phi )^{2r + s }.
\eeq
So we may assume that each monomial in (\ref{eqn:sbd}) is maximally expanded.  As in \cite{Kevin1} we observe that for every single label $a \in S ( \Gamma )$ there is a label $b \in \{1, ..., 2r \} \backslash \{ a \}$ s.t. the monomial $ ( F_{i_b} )_{\sigma_b}$ contains an off-diagonal resolvent entry with $i_a$ as a lower index.

Hence, 
\beq
\sum_{k=1}^{2r} b ( \sigma_k ) \geq s
\eeq
and we get the claim from (\ref{eqn:onesbd}).

For our fixed partition $\Gamma$ denote its size by $l = | \Gamma |$.  It follows from (\ref{eqn:sbd}) that
\beq
\left| \frac{1}{N^{2r} } \sum_{i_1, ..., i_r} \1_{ \{ \Gamma = \Gamma ( \underline{i} ) \} } \ee [ g_{i_1}^2 Q_{i_1} ( \widebar{ G_{i_1, i_1} ^{-1} } ) ... g_{i_r}^2 Q_{i_r} (  G_{i_r, i_r} ^{-1}  ) ]  \right| \leq C \frac{1}{N^{2r}} \sum_{k_1=1}^N ... \sum_{k_l =1}^N |g_{k_1} |^{d_1}... |g_{k_l } |^{d_l} ( \varphi^{2 \xi} \Phi )^{ 2r+s}
\eeq
where $d_i$ is the size of the $i$th block of the partition $\Gamma$.  Using $|g_k| \leq \eta^{-1}$ we then bound the sum
\begin{align}
\frac{1}{N^{2r} }\sum_{k_1=1}^N ...  \sum_{k_l =1}^N |g_{k_1} |^{d_1}... |g_{k_l } |^{d_l} &\leq \left( \frac{1}{ N \eta } \right)^{2r - l } \frac{1}{N^l}\sum_{k_1=1}^N ...  \sum_{k_l =1}^N |g_{k_1} |... |g_{k_l } | \notag \\
&\leq \left( \frac{1}{ N \eta } \right)^{2r - l }  \log (N)^{l} \leq  \varphi^{ \xi} (\Phi )^{ 4r - 2 l }.
\end{align}
We have the inequality $2r + s +  (4r - 2 l)  \geq 4 r$.  Inserting the two bounds just derived into the RHS of (\ref{eqn:bigsum}) we get
\beq
\ee \left| \frac{1}{N} \sum_{i=1}^N g_i^2 Q_i \left( \frac{1}{ G_{ii} } \right) \right|^{2r} \leq \sum_{ \Gamma \in \P_{2r} } ( \varphi^{ 3\xi } \Phi )^{4 r } \leq (C r)^{C r }  ( \varphi^{ 3 \xi } \Phi )^{4 r }.
\eeq
\qed

\subsubsection{Proof of local law in $\D_2$}

The proof of the local law in $\D_2$ is similar but shorter than the proof in $\D_1$.  The starting point is once again the Schur complement formula which yields
\begin{align}
m_N - \mfc &= \frac{1}{N} \sum_{i=1}^N \frac{1}{ V_i  - z - T \mfc (z) + T \psi + Y_i } - \frac{1}{ V_i  - z - T \mfc (z) } \notag \\
&= \frac{1}{N} \sum_{i=1}^N \frac{ - T \psi - T Y_i }{ (V_i - z - T \mfc (z) + T \psi + Y_i )  (V_i - z - T \mfc (z) )}.
\end{align}
From the fact that $\eta \geq 10$ on $\D_2$ and that $|Y_i| \leq \varphi^\xi / N^{1/2}$ with $(\xi, \nu)$-high probability we obtain
\beq
|m_N - \mfc | \leq \frac{ | m_N - \mfc |}{2} + \varphi^{\xi} \frac{1}{ \sqrt{N} }
\eeq
with $(\xi, \nu)$-high probability and so $| m_N - \mfc | \leq \varphi^{\xi} /\sqrt{N}$ with $(\xi, \nu)$-high probability, uniformly in $\D_2$.

We now proceed as in the proof of Lemma \ref{lem:sself} and expand the Schur complement formula and obtain
\begin{align}
m_N &= \frac{1}{N} \sum_{i=1}^N g_i + \frac{1}{N} \sum_{i=1}^N g_i^2 ( T \psi + TY_i) + \O \left(  \frac{1}{N} \sum_{i=1}^N  |g_i|^3 \varphi^{\xi} \frac{1}{N} \right)
\notag \\
&= \frac{1}{N} \sum_{i=1}^N g_i + \frac{1}{N} \sum_{i=1}^N g_i^2 ( T \psi + Q_i [ G_{ii}^{-1} ]) + \O \left(  \frac{1}{N} \sum_{i=1}^N  |g_i|^3 \varphi^{\xi} \frac{1}{N} \right).
\end{align}
In the first line we used that $|\psi| + |Y_i| \leq \varphi^{\xi}/\sqrt{N}$ with $(\xi, \nu)$-high probability, and in the second line that $Y_i = Q_i [ G_{ii}^{-1} ] + \O ( (N \eta)^{-1} )$.  Since $| 1 - T R_2 | \geq c$ uniformly in $\D_2$ this yields with $(\xi, \nu)$-high probability
\beq
|m_N - \mfc | \leq  \frac{g(z)  \varphi^\xi }{N} + C \left| \frac{1}{N} \sum_{i=1}^N g_i^2 Q_i [ G_{ii}^{-1} ]\right|
\eeq
where we defined
\beq
g(z) := \frac{1}{N} \sum_{i=1}^N \frac{ 1}{ (V_i - E )^2 + \eta^2 } = \frac{1}{\eta} \Im [ m_V (z) ].
\eeq
Since $g(z) \leq \eta^{-2}$ we see it suffices to show that
\beq \label{eqn:dtwofluct}
 \left| \frac{1}{N} \sum_{i=1}^N g_i^2 Q_i [ G_{ii}^{-1} ] \right| \leq \varphi^{\xi} \frac{1}{ N \eta^2}
\eeq
with $(\xi, \nu)$-high probability uniformly in $\D_2$.  Before proving \eqref{eqn:dtwofluct} we note that the proof of Lemma \ref{lem:apriori2} yields
\bel
Uniformly for $| \tt | \leq ( \log N)^{\xi-1}$ and $z \in \D_2$ we have
\beq \label{eqn:aptwo1}
\left| \frac{g_i}{ G_{ii}^{\btt} } - 1 \right|  \leq \frac{ \varphi^{\xi} }{ \sqrt{N}},
\eeq
and for $i \neq j$
\beq \label{eqn:aptwo2}
\left| G_{ij}^{\btt} \right| \leq |g_i g_j | \frac{ \varphi^\xi}{ \sqrt{N}}
\eeq
and
\beq \label{eqn:aptwo3}
\left| Q_i \left[ \frac{1}{ G_{ii}^{\btt} } \right] \right| \leq \frac{ \varphi^\xi} { \sqrt{N}}.
\eeq
\eel

\noindent{\bf Proof of \eqref{eqn:dtwofluct}}.  We proceed as in the proof of Lemma \ref{lem:fa}.  We use the same notation appearing there  and will not redefine it.  As in the proof of Lemma \ref{lem:fa}, we start by writing
\beq  \label{eqn:bsumtwo}
\ee \left| \frac{1}{N} \sum_{i=1}^N g_i^2 Q_i \left( \frac{1}{ G_{ii} } \right) \right|^{2r} = \frac{1}{N^{2r} } \sum_{ \Gamma \in \P_{2r} } \sum_{i_1, ..., i_r} \1_{ \{ \Gamma = \Gamma ( \underline{i} ) \} } \ee [ g^2_{i_1} Q_{i_1} ( \widebar{ G_{i_1, i_1} ^{-1} } ) ... g^2_{i_r} Q_{i_{2r}} (  G_{i_r, i_r} ^{-1}  ) ].
\eeq
We apply the same expansion algorithm as in the proof of Lemma \ref{lem:fa} to each term in the above sum and estimate the resulting terms. With $b ( \sigma_k)$ denoting the number of ones in the string $\sigma_k$, we claim that
\beq
| Q_{i_k} ( F_{i_k} )_{\sigma_k } | \leq ( \varphi^\xi N^{-1/2} )^{1 + b ( \sigma_k ) }.
\eeq
When $b ( \sigma_k ) = 0$ this follows from \eqref{eqn:aptwo3}.  For the moment let $\ell ( \sigma_k)$ denote the number of off-diagonal resolvent entries generated by the algorithm in the numerator.  If $b( \sigma_k )>0$, note that $\ell ( \sigma_k ) \geq b ( \sigma_k ) + 1$.  Note also that there are $\ell ( \sigma_k ) + 1$ diagonal resolvent entries in the denominator.  We want to use the bound $\eqref{eqn:aptwo2}$ and the bound $| G_{ii}^{\btt} |^{-1} \leq 2 | g_i |^{-1}$ to conclude
\beq \label{eqn:onesbd2}
| Q_{i_k} ( F_{i_k} )_{\sigma_k } | \leq ( \varphi^\xi N^{-1/2} )^{\ell ( \sigma_k ) }.
\eeq
The problem is that $|g_i|$ can be quite small and so our estimate $|G_{ii}|^{-1} \leq 2 |g_{i}|^{-1}$ can be quite big, and we need to compensate for this with the factor $|g_i g_j|$ that appears in the estimate \eqref{eqn:aptwo2} for the off-diagonal entries.

For the monomial $ ( F_{i_k })_{\sigma_k }$ define $p_j ( (F_{i_k } )_{\sigma_k } )$ to be the number of times a diagonal resolvent entry with lower index $j$ appears in $(F_{i_k} )_{\sigma_k}$.  Define $s_j ( ( F_{i_k })_{\sigma_k } )$ to be the number of times that $j$ appears as a lower index of any off-diagonal resolvent entry.  For example, for the monomial
\beq
\frac{ G_{ik}^\btt G_{ki}^\btt}{ G_{ii}^\btt G_{ii}^\bttk G_{kk}^\btt }
\eeq
we have $s_i = 2$, $s_k = 2$, $p_i = 2$ and $p_k = 1$.  For \eqref{eqn:onesbd2} to hold, we need $p_i \leq s_i$ for every $i$ and every monomial.  It is easy to see that by definition of the expansion algorithm that we will always have $p_i \leq s_i$ at each stage of the algorithm as soon as a $1$ appears in $\sigma_k$.  Therefore, \eqref{eqn:onesbd2} holds.

Estimating $|g_i | \leq \eta^{-1}$ we see that, as in the proof of Lemma \ref{lem:fa} that any nonzero term on the RHS of \eqref{eqn:bsumtwo} satisfies
\beq
| g_{i_1} \widebar{ Q_{i_1} ( F_{i_1} )_{\sigma_1} }\cdots g_{i_{2r}} Q_{i_{2r}} ( F_{i_{2r} } )_{\sigma_{2 r } } | \leq \eta^{- 4 r} ( \varphi^{2 \xi } N^{-1/2} )^{2 r + s},
\eeq
where $s$ is the number of single labels appearing in the term.  It then follows that, with $l = | \Gamma |$,
\begin{align}
\left| \frac{1}{N^{2r} } \sum_{i_1, ..., i_r} \1_{ \{ \Gamma = \Gamma ( \underline{i} ) \} } \ee [ g_{i_1}^2 Q_{i_1} ( \widebar{ G_{i_1, i_1} ^{-1} } ) ... g_{i_r}^2 Q_{i_r} (  G_{i_r, i_r} ^{-1}  ) ]  \right| &\leq C \frac{\eta^{- 4 r}}{N^{2r}} \sum_{k_1=1}^N ... \sum_{k_l =1}^N ( \varphi^{2 \xi} N^{-1/2} )^{ 2r+s} \notag \\
&\leq C \frac{ \eta^{- 4 r}}{ N^{2r} } \varphi^{ 4 r \xi } \frac{ N^l}{ N^{r+s/2} } \notag \\
&\leq C \eta^{-4 r } \varphi^{4 r \xi } \frac{1}{ N^{2 r }}
\end{align}
where we used that $s/2 + r - l \geq 0$.  As in the proof of Lemma \ref{lem:fa} we now get the claim. \qed

\subsection{Proof of Theorem \ref{thm:rig}}
\label{sect:rigsketch}

In this section we show how the Helffer-Sj\"ostrand formula is used to yield rigidity for the eigenvalues of $H$.  For the most part the material in this section is standard; see, e.g., \cite{general}.  In this section we will always work on the event that the statement of Theorem \ref{thm:ll} holds.  Let $f$ be a smooth compactly supported function.  We have
\beq \label{eqn:hs}
f( \lambda ) = \frac{1}{ \pi } \int_{\rr^2} \frac{ \i y f'' (x) \chi (y) + \i ( f(x) + \i y f' (x) ) \chi' (y) }{ \lambda - x - \i y } \d x \d y
\eeq
where $\chi$ is a smooth compactly supported function that is $1$ in a neighborhood of the origin.  Fix now $E_2 \in \I_{E_0, q G}$ and let $E_1 := - N^{4 B_V}$ (recall the definition of $B_V$ in Definition \ref{def:v}).  Take $f$ to satisfy
\beq
f (E) = 0, \quad E \notin [ E_1 - 1, E_2 +  \eta ], \qquad f(E) = 1, \quad E \in [ E_1, E_2] ,
\eeq
where
\beq
\eta = \frac{ \varphi^{ 50 \xi}}{N}.
\eeq
We can take $|f' (E) |, |f'' (E) | \leq C$ for $E$ near $E_1$ and $ | f' (E) | \leq C \eta$ and $| f''(E)  | \leq C \eta^{-2}$ for $E$ near $E_2$. For this smoothed out eigenvalue counting function we prove the following estimate.
\bel \label{lem:smoothed} Let $f$ be as above.  There exists a $c_2 >0$ such that with $(\xi, \nu)$-high probability we have
\beq
\left| \int f ( \lambda ) \d \rho_N ( \lambda ) - \int f ( \lambda ) \d \rhofc ( \lambda ) \right| \leq \frac{ \varphi^{ c_2 \xi } }{N},
\eeq
where $\rho_N$ denotes the empirical measure of $H_T$.
\eel
\proof 
Take $\chi$ a smooth cut-off function to satisfy
\beq
\chi (y) = 1, \quad |y| \leq N^{10 B_V}, \qquad \chi (y) = 0, \quad |y| > N^{10 B_V} + 1.
\eeq
It is no loss of generality to assume that $B_V \geq 10$. 
Let $S = m_N - \mfc$ and $\tilrho = \rho_N - \rhofc$.  From \eqref{eqn:hs} we have
\begin{align} \label{eqn:hserror}
\left| \int f( \lambda ) \tilrho ( \lambda ) \d \lambda \right| &= \left| \Re \int f( \lambda ) \tilrho ( \lambda ) \d \lambda \right| \notag \\
&\leq C \left| \int y f'' (x) \chi (y) \Im S (x + \i y ) \d x \d y \right| \notag \\
&+ C \int | f (x) \chi' (y) | | \Im S (x + \i y ) | + |y f'(x) \chi' (y) | | \Re S ( x + \i y ) | \d x \d y.
\end{align}
Since $| S (x + \i y ) | \leq \varphi^{c_1 \xi } / (N y)$ for $y \geq 10$ we immediately get
\beq
\int |y| | f' (x) | | \chi' (y) | \Re S (x + \i y ) | \d x \d y \leq C \frac{  \varphi^{ c_1 \xi } }{N},
\eeq
where we used that $| f' (x) | \leq \eta^{-1}$ is nonzero only on an interval of length $\eta$ near $E_2$, and is bounded near $E_1$. 
Since $\chi' (y) \neq 0$ only for $N^{10 B_V } \leq |y| \leq N^{ 10 B_V} +1$ and $f(x) \neq 0$ only on an interval of length $2 N^{4 B_V}$ we have
\beq
\int |f (x) \chi' (y) | | \Im S (x + \i y ) | \d x \d y \leq \frac{ N^{5 B_V}}{ N^{ 10 B_V } } \leq \frac{1}{N}.
\eeq

We bound the second line of \eqref{eqn:hserror} as
\begin{align}
 \left| \int y f'' (x) \chi (y) \Im S (x + \i y ) \d x \d y \right| &\leq  \left| \int_{ |x - E_1 | \leq 2, |y| \leq 10 } y f'' (x) \chi (y) \Im S (x + \i y ) \d x \d y \right| \notag \\
 &+  \left| \int_{ |x-E_1| \leq 2, |y| > 10} y f'' (x) \chi (y) \Im S (x + \i y ) \d x \d y \right| \notag \\
&+  \left| \int_{ |x - E_2| \leq 2, |y| \leq \eta} y f'' (x) \chi (y) \Im S (x + \i y ) \d x \d y \right| \notag \\
&+  \left| \int_{ |x - E_2 | \leq 2, |y| > \eta } y f'' (x) \chi (y) \Im S (x + \i y ) \d x \d y \right| \notag \\
&=: A_1 + A_2 + A_3 +A_4
\end{align}
We easily see that
\beq
|A_1| \leq N^{ - B_V}
\eeq
using that $\dist (E_1, \supp \tilrho ) \geq N^{3 B_V }$.  For $A_2$ we integrate by parts in $x$ and then use $\del_x \Im S = - \del_y \Re S$ to integrate by parts in $y$ to obtain
\begin{align}
A_2 &\leq \left| \int_{ |x - E_1 | \leq 2 ,  |y| > 10} f' (x) \del_y ( y \chi (y) ) \Re S (x + \i y ) \d x \d y \right| \notag \\
&+ \left| \int_{ |x - E_1 | \leq 2 } f' (x) 10 \chi (10) \Re S ( x + 10 \i  ) \d x \right| =: B_1 + B_2.
\end{align}
Using $|S| \leq \varphi^{c_1 \xi}/ ( N y)$ we see the first term is bounded by
\beq
B_1 \leq C \frac{  \log (N) \varphi^{ c_1 \xi } }{N}.
\eeq
The term $B_2$ is easily bounded by $C \varphi^{ c_1 \xi } /N$.

For the term $A_3$, note that since that $y \to y \Im [ m (E + \i y )]$ is a decreasing function for any Stieltjes transform $m$ of a positive measure, we see that, using that $| \mfc|$ is bounded near $E_2$,
\beq
\Im [ m_N ( x + \i y ) ] \leq \frac{ \eta}{y} \Im [ m_N (x + \i \eta ) ] \leq C \frac{ \varphi^{ 50} \xi }{ N y }, \qquad y < \eta
\eeq
and so $\Im [ S(x + \i y ) ] \leq C \varphi^{ c_1 \xi } / (N y )$ for $y < \eta$.  Using then that the $x$ integration is over an interval of length $\eta$ and $|f'' | \leq C \eta^{-2}$ we see that
\beq
A_3 \leq C \eta \leq C \frac{ \varphi^{ 50 \xi } }{ N }.
\eeq
For $A_4$ we integrate by parts as before and obtain
\begin{align}
A_4 &\leq \left| \int_{|y| > \eta , |x - E_2| \leq 2 } f' (x) \del_y ( y \chi (y) ) \Re S (x + \i y ) \d x \d y \right| \notag \\
&+ \left| \int_{ |x - E_2 | \leq 2 } f' (x) \eta \chi ( \eta ) \Re S (x + \i y ) \d x \right| =: B_3 + B_4.
\end{align}
The term $B_3$ is estimated using $|S| \leq \varphi^{ c_1 \xi}/ (Ny )$.  For the term $B_4$ we use that the $x$ integration is over an interval of length $\eta$, that $|f'(x) | \leq C \eta^{-1}$ and $|S| \leq \varphi^{c_1 \xi} / (N \eta )$ to obtain
\beq
B_4 \leq C \frac{ \varphi^{ c_1 \xi }}{ N }  .
\eeq
This yields the claim. \qed

Define the empirical eigenvalue counting function and the eigenvalue counting function of the deformed semicircle law by
\beq
\nt (E) = \frac{1}{N} \left| \{ i : \lambda_i ( \Ht ) \leq E \} \right| , \qquad \nfct (E) = \int_{ -\infty}^E \rhofct (E' ) \d E' 
 \eeq
 respectively.  
 From Lemma \ref{lem:smoothed} we conclude the following lemma.
 \bel \label{lem:counting} Let $q> 0$ be as above and let $\om > 0$ and $t \in \Tom$.  With $(\xi, \nu)$-high probability we have uniformly for $E \in \I_{E_0, qG}$,
 \beq
 \left| \nt (E) - \nfct (E) \right| \leq \frac{ \varphi^{ c_3 \xi } } { N}
 \eeq
 for a constant $c_3 > 0$ depending on $\IV$, $A_0$, $C_1$ in (\ref{eqn:param}), the constants appearing in (\ref{eqn:vass}) and the choice of $q$. 
 \eel
 \proof Let $f$ be as in Lemma \ref{lem:smoothed}.  Note that spectrum of $H_T$ is contained in $[-N^{ B_V} -1, N^{B_V} +1 ]$ with $(\xi, \nu)$-high probability.
 The left endpoint $E_-$ of the support of $\rhofct$ is the smallest solution to the equation
 \begin{align}
 T^{-1} = \frac{1}{N} \sum_{i=1}^N \frac{1}{ (V_i - E_- - T \mfcll ( E ) )^2 }.
 \end{align}
 Clearly $E_- \geq -N^{B_V}-1$. With $(\xi, \nu)$-high probability we then have that, with $\eta$ as in Lemma \ref{lem:smoothed},
  \begin{align}
 \nt (E_2) &\leq \int f (E') \d \nt (E' ) \leq \int  f (E' ) \d \rhofct (E' ) +  \frac{ \varphi^{c_2 \xi }}{N} \notag \\
 & \leq \nfct (E_2 + \eta ) +  \frac{ \varphi^{c_2 \xi }}{N} \leq \nfct (E_2) + \frac{ \varphi^{c_2 \xi }}{N}  + C \eta \notag \\
 &\leq \nfct (E_2) + \frac{ \varphi^{ c_3 \xi }}{N}
 \end{align}
 where we used the fact that $\rhofct \leq C$ in $\I_{E_0, q G}$.  The lower bound is similar. \qed

The conclusion of Theorem \ref{thm:rig} from Lemma \ref{lem:counting} is standard and is similar to, for example, the proof of Theorem 7.6 in \cite{general}.

\bibliography{mybib}{}
\bibliographystyle{abbrv}

\end{document}